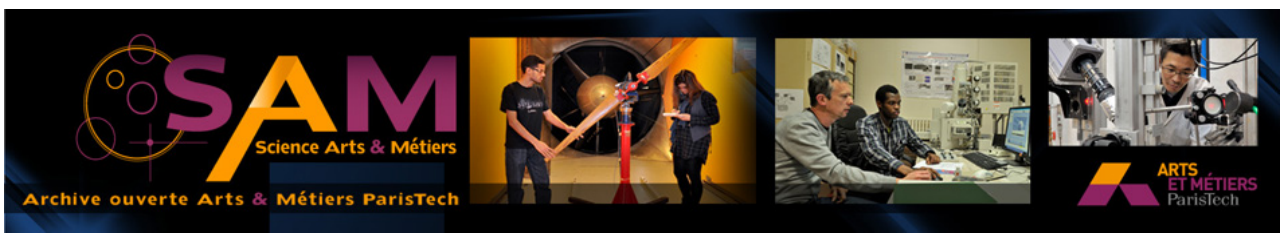

# A simple Cartesian scheme for compressible multimaterials

Yannick Gorsse [a,b,c], Angelo Iollo [a,b,c,*], Thomas Milcent [d,e], Haysam Telib [f]

### ABSTRACT

We present a simple numerical method to simulate the interaction of two non-miscible compressible materials separated by an interface. The media considered may have significantly different physical properties and constitutive laws, describing for example fluids or hyperelastic solids. The model is fully Eulerian and the scheme is the same for all materials. We show stiff numerical illustrations in case of gas–gas, gas–water, gas–elastic solid interactions in the large deformation regime.

## 1. Introduction

Physical phenomena that involve several materials are ubiquitous in nature and applications: multiphase flows, fluid–structure interaction, impacts, to cite just a few examples. In recent literature several strategies were proposed to attack these problems: Lagrangian models [29], Arbitrary-Lagrangian–Eulerian (ALE) models [11,20], Eulerian models [4,9,19,23]. In this context, immersed boundary methods [17,21] are an option to discretize boundary conditions at the interface, trading off accuracy for mesh generation simplicity. The literature in this domain is vast and each of the approaches cited include a variety of methods that are adapted to specific applications.

Broadly speaking Lagrangian methods are interesting because the interface between the materials is fixed in the reference domain and a body-fitted mesh is generated once for all. The interface equilibrium conditions can be taken into account in a simple and accurate way. Free-surface problems and purely elastic problems are the paradigmatic applications of this approach. On the other hand, if one of the materials is a fluid, the mapping from the physical to the reference domain can become very irregular or even singular for large times. ALE methods take into account this peculiarity of fluid flows by discretizing the physical domain over an unsteady mesh that can also cope with the moving interfaces. Moreover, the mapping from the physical domain to the computational domain is independent of the trajectories of the fluid particles and therefore its regularity can be kept under control. The counterpart is a more complex scheme formulation and numerical implementation. Also, when large deformations occur, mesh generation and partitioning can pose challenging problems. In Eulerian methods the problem is approximated on a fixed mesh in the physical domain, but additional modeling is required to obtain a consistent, stable and accurate description of the (eventually smoothed) material interface evolution.

---

\* Corresponding author at: Univ. Bordeaux, IMB, UMR 5251, F-33400 Talence, France.
  *E-mail addresses:* yannick.gorsse@math.u-bordeaux1.fr (Y. Gorsse), angelo.iollo@math.u-bordeaux1.fr (A. Iollo), thomas.milcent@u-bordeaux1.fr

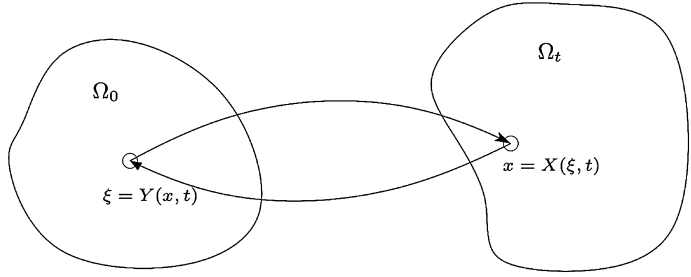

**Fig. 1.** Forward and backward characteristics.

Immersed boundary methods can be used to keep the material interface sharp. With this approach, the contact discontinuity can arbitrarily cross the grid and the transmission conditions, usually in terms of velocity and stress, are applied via interpolation. The models considered on either side of the interface can be either Lagrangian or Eulerian, using Lagrangian markers or level-set functions to describe the interface according to the specific problem considered.

In this study, we are interested in the numerical simulation of phenomena such as the scattering of shock waves at gas–water interfaces, the propagation of non-linear elastic waves from a hyperelastic solid to a fluid and vice versa, the scattering of non-linear elastic waves at solid–solid material discontinuities. These phenomena can be modeled by a fully Eulerian system of conservation laws that applies to every material; only the constitutive law may change, reproducing the mechanical characteristics of the medium under consideration. For example, an elastic material or a gas will be modeled by the same set of quasi-linear hyperbolic partial differential equations (PDEs) except for the constitutive law relating the material deformation and the stress tensor.

The systematic derivation of such models starting from continuum mechanics principles, their thermodynamic consistency and the corresponding wave-propagation patterns were initially studied in [13]. Their numerical simulation is delicate because standard Godunov schemes lead to pressure oscillations at the material contact discontinuity already in the case of multifluids. In [1] the pressure perturbation mechanism at the origin of this phenomenon was explained and a first fix was proposed. An effective remedy to this problem was presented in [10] with the ghost-fluid method (sharp interface between the materials). For multifluids, improvements of this approach requiring less storage were proposed in [2] (diffuse interface) and [8] (sharp interface). The common idea of these methods is to define a "ghost" fluid that has continuous mechanical characteristics across the interface, but the same thermodynamic state or the same equation of state of the actual fluid. This assumption leads to locally non-conservative schemes that are consistent, stable and non-oscillatory at the material interface.

For elastic compressible materials existing methods either rely on the definition of ghost materials (see [4,23] for hyperelastic models and [19,27] for a hypoelastic formulation) or on mixture models and diffuse interfaces [9]. In this paper we propose a simple, stable and non-oscillatory scheme for multimaterials that avoids the definition of a ghost hyperelastic medium. The equilibrium boundary conditions at the material interface are imposed like in immersed boundary methods, in the same spirit of what is done in [14] for rigid bodies. Therefore, in our scheme the evolution of the material discontinuity is sharp by construction. In this paper we limit the discussion to the two dimensional case but the models and the numerical scheme can be extended in 3D.

In the following we detail the Eulerian model, the numerical method and we present a set of stiff test-cases involving fluid and hyperelastic compressible materials.

## 2. The model

This model was already discussed in [6,9,13,23–25]. In this section, we develop the principal elements of the formulation.

### 2.1. Forward and backward characteristics

Let $\Omega_0 \subset \mathbb{R}^2$ be the reference or initial configuration of a continuous medium and $\Omega_t \subset \mathbb{R}^2$ the deformed configuration at time $t$. In order to describe the evolution of this medium in the Lagrangian frame we define the forward characteristics $X(\xi, t)$ as the image at time $t$ in the deformed configuration of a material point $\xi$ belonging to the initial configuration, i.e., $X : \Omega_0 \times [0, T] \to \Omega_t$, $(\xi, t) \mapsto X(\xi, t)$ (see Fig. 1). The corresponding Eulerian velocity field $u$ is defined as $u : \Omega_t \times [0, T] \to \mathbb{R}^2$, $(x, t) \mapsto u(x, t)$ where

$$\begin{cases} X_t(\xi, t) = u\big(X(\xi, t), t\big) \\ X(\xi, 0) = \xi \end{cases} \tag{1}$$

To describe the continuous medium in the Eulerian frame, we introduce the backward characteristics $Y(x, t)$ (see [6]) that for a time $t$ and a point $x$ in the deformed configuration, gives the corresponding initial point $\xi$ in the initial configuration, i.e., $Y : \Omega_t \times [0, T] \to \Omega_0$, $(x, t) \mapsto Y(x, t)$ (see Fig. 1). Since $Y(X(\xi, t), t) = \xi$, differentiating with respect to time and space in turn we have:

$$\begin{cases} Y_t + u \cdot \nabla_x Y = 0 \\ Y(x, 0) = x \end{cases} \tag{2}$$

and

$$\left[ \nabla_\xi X(\xi, t) \right] = \left[ \nabla_x Y(x, t) \right]^{-1}. \tag{3}$$

The relation (2) is the Eulerian equivalent of the characteristic equation (1). In addition, Eq. (3) allows to compute the gradient of the deformation in the Eulerian frame via $Y$. We are now in the position of writing the mass, momentum and energy conservation with respect to both the deformed (Eulerian) and the reference (Lagrangian) configuration.

## 2.2. Eulerian and Lagrangian formulations

The integral form of the governing equations of mass, momentum and energy conservation in the deformed configuration $\Omega_t$ can be written as

$$\begin{cases} \dfrac{\mathrm{d}}{\mathrm{d}t} \left( \displaystyle\int_{\Omega_t} \rho \, \mathrm{d}x \right) = 0 & \text{conservation of mass} \\[2ex] \dfrac{\mathrm{d}}{\mathrm{d}t} \left( \displaystyle\int_{\Omega_t} \rho u \, \mathrm{d}x \right) = \displaystyle\int_{\partial\Omega_t} \sigma n \, \mathrm{d}s & \text{conservation of momentum} \\[2ex] \dfrac{\mathrm{d}}{\mathrm{d}t} \left( \displaystyle\int_{\Omega_t} \rho e \, \mathrm{d}x \right) = \displaystyle\int_{\partial\Omega_t} \left( \sigma^T u \right) \cdot n \, \mathrm{d}s & \text{conservation of energy} \end{cases} \tag{4}$$

where $\partial\Omega_t$ is the boundary of $\Omega_t$ and $n(x, t)$ is the corresponding outward normal unitary vector. The physical variables are the density $\rho(x, t)$, the velocity $u(x, t)$, the total energy per unit mass $e(x, t)$ and the Cauchy stress tensor in the physical domain $\sigma(x, t)$. With the change of variables $x = X(\xi, t)$, the mass equation becomes in the reference configuration $\Omega_0$

$$\rho \left( X(\xi, t), t \right) \det \left( \nabla_\xi X(\xi, t) \right) = \rho_0(\xi) \tag{5}$$

where $\rho_0(\xi)$ is the initial density. The counterpart of momentum and energy equations can be written in the reference configuration $\Omega_0$ with the change of variables $x = X(\xi, t)$ and (5) as

$$\begin{cases} \dfrac{\mathrm{d}}{\mathrm{d}t} \left( \displaystyle\int_{\Omega_0} \rho_0 u \, \mathrm{d}\xi \right) = \displaystyle\int_{\partial\Omega_0} \mathcal{T} n_0 \, \mathrm{d}s_0 \\[2ex] \dfrac{\mathrm{d}}{\mathrm{d}t} \left( \displaystyle\int_{\Omega_0} \rho_0 e \, \mathrm{d}\xi \right) = \displaystyle\int_{\partial\Omega_0} \left( \mathcal{T}^T u \right) \cdot n_0 \, \mathrm{d}s_0 \end{cases} \tag{6}$$

where $\partial\Omega_0$ is the boundary of $\Omega_0$ and $n_0(\xi)$ is the corresponding outward normal unitary vector. Here $\mathcal{T}(\xi, t)$ is the first Piola–Kirchhoff stress tensor and is related to the Cauchy stress tensor $\sigma(x, t)$ by the Piola transformation

$$\mathcal{T}(\xi, t) = \sigma \left( X(\xi, t), t \right) \mathrm{Cof} \left( \nabla_\xi X(\xi, t) \right) \tag{7}$$

where the cofactor matrix is defined by $\mathrm{Cof}(A) = \det(A) A^{-T}$. The local form of Eqs. (4) in the deformed configuration is

$$\begin{cases} \rho_t + \mathrm{div}_x(\rho u) = 0 \\ (\rho u)_t + \mathrm{div}_x(\rho u \otimes u - \sigma) = 0 \\ (\rho e)_t + \mathrm{div}_x(\rho e u - \sigma^T u) = 0 \end{cases} \tag{8}$$

The local form of Eqs. (6) in the reference configuration is

$$\begin{cases} (\rho_0 u)_t + \mathrm{div}_\xi(-\mathcal{T}) = 0 \\ (\rho_0 e)_t + \mathrm{div}_\xi(-\mathcal{T}^T u) = 0 \end{cases} \tag{9}$$

To close these systems a constitutive law is needed to detail the relationship between the stress tensors and the other physical variables.

## 2.3. Hyperelastic models

The classical results of this section can be found in the textbook [16]. Let us define $\varepsilon = e - \frac{1}{2}|u|^2$ the internal energy per unit mass. In the hyperelastic context, $\varepsilon$ is a function of the strain tensor $\nabla_\xi X$ and the entropy $s$, so that

$$\varepsilon = \varepsilon\big(\nabla_\xi X(\xi,t), s\big(X(\xi,t),t\big)\big) \tag{10}$$

By definition, the Piola–Kirchhoff stress tensor can be written as the derivative of $\varepsilon(F,s)$ with respect to the first variable at fixed entropy, i.e.,

$$\mathcal{T}(\xi,t) = \rho_0 \frac{\partial \varepsilon}{\partial F}\bigg|_s (\nabla_\xi X, s) \tag{11}$$

The energy has to be Galilean invariant and we focus in this paper on the isotropic case. It can be proved that the material is Galilean invariant and isotropic if, and only if, $\varepsilon$ is expressed as a function of the invariants of the right Cauchy–Green tensor $C(\xi,t) = [\nabla_\xi X]^T [\nabla_\xi X]$ or equivalently of the invariants of the left Cauchy–Green tensor $B(\xi,t) = [\nabla_\xi X][\nabla_\xi X]^T$. The two dimensional invariants often considered in the literature are $\mathrm{Tr}(\cdot)$ and $\mathrm{Det}(\cdot)$.

We assume that $\varepsilon$ is the sum of $\varepsilon_{\mathrm{vol}}$, a term depending on volume variation and entropy, and $\varepsilon_{\mathrm{iso}}$ a term accounting for isochoric deformation. In general the term relative to an isochoric transformation may also depend on entropy. Here, we will limit the discussion to materials where the isochoric term is independent of the entropy, as in the case of idealized crystals (metals, ceramic).

The internal energy is given by

$$\varepsilon = \varepsilon_{\mathrm{vol}}\big(\rho\big(X(\xi,t),t\big), s\big(X(\xi,t),t\big)\big) + \varepsilon_{\mathrm{iso}}\big(\mathrm{Tr}\big(\overline{B}(\xi,t)\big)\big) \tag{12}$$

With the change of variables $\xi = Y(x,t)$ and (3) we obtain the Eulerian equivalent of (12)

$$\varepsilon = \varepsilon_{\mathrm{vol}}\big(\rho(x,t), s(x,t)\big) + \varepsilon_{\mathrm{iso}}\big(\mathrm{Tr}\big(\overline{B}(x,t)\big)\big) \tag{13}$$

where in 2D

$$\overline{B}(x,t) = [\nabla_x Y]^{-1} [\nabla_x Y]^{-T} / J(x,t) \qquad J(x,t) = \det\big([\nabla_x Y]\big)^{-1} \tag{14}$$

It can be remarked that $\overline{B}$ accounts for isochoric deformations since $\det(\overline{B}) = 1$. The formula (11) allows to compute $\mathcal{T}(\xi,t)$ in the case of the constitutive law (12) with (5) and we deduce $\sigma(x,t)$ with (7) and the change of variables $\xi = Y(x,t)$ (see [16] for the details). The result is

$$\sigma(x,t) = -\rho^2 \frac{\partial \varepsilon_{\mathrm{vol}}}{\partial \rho}\bigg|_s (\rho,s)I + 2J^{-1}\varepsilon'_{\mathrm{iso}}\big(\mathrm{Tr}(\overline{B})\big)\left(\overline{B} - \frac{\mathrm{Tr}(\overline{B})}{2}I\right) \tag{15}$$

## 2.4. Summary of the hyperelastic Eulerian model used

Together with equations of mass, momentum and energy conservation (8), the additional equation (2) is required in order to record the deformation in the Eulerian frame. However, since $\sigma$ will directly depend on $\nabla_x Y$ (13)–(15) we take the gradient of (2) as a governing equation and obtain the system in conservative form

$$\begin{cases} \rho_t + \mathrm{div}_x(\rho u) = 0 \\ (\rho u)_t + \mathrm{div}_x(\rho u \otimes u - \sigma) = 0 \\ (\nabla_x Y)_t + \nabla_x(u \cdot \nabla_x Y) = 0 \\ (\rho e)_t + \mathrm{div}_x\big(\rho e u - \sigma^T u\big) = 0 \end{cases} \tag{16}$$

where some fluxes of $\nabla_x Y$ are 0 (see Eq. (20)). The initial density $\rho(x,0)$, the initial velocity $u(x,0)$, the initial total energy $e(x,0)$ and $\nabla_x Y(x,0) = I$ are given together with appropriate boundary conditions.

We choose a general constitutive law that models gas, fluids and elastic solids. The internal energy per unit mass $\varepsilon = e - \frac{1}{2}|u|^2$ is defined as

$$\varepsilon(\rho, s, \nabla_x Y) = \underbrace{\frac{\kappa(s)}{\gamma - 1}\left(\frac{1}{\rho} - b\right)^{1-\gamma} - a\rho + \frac{p_\infty}{\rho}}_{\text{general gas}} + \underbrace{\frac{\chi}{\rho_0}\big(\mathrm{Tr}(\overline{B}) - 2\big)}_{\text{neohookean solid}} \tag{17}$$

where $\overline{B}$ is defined in (14). We obtain according to (15)

$$\sigma(\rho, s, \nabla_x Y) = -p(\rho,s)I + 2\chi J^{-1}\left(\overline{B} - \frac{\mathrm{Tr}(\overline{B})}{2}I\right) \tag{18}$$

**Table 1**
Typical material parameters.

| Material | $\gamma$ | $a$ [Pa kg$^{-2}$ m$^6$] | $b$ [kg$^{-1}$ m$^3$] | $p_\infty$ [Pa] | $\chi$ [Pa] |
|---|---|---|---|---|---|
| Perfect gas (Air) | 1.4 | 0 | 0 | 0 | 0 |
| Van Der Waals gas | 1.4 | 5 | $10^{-3}$ | 0 | 0 |
| Stiffened gas (Water) | 4.4 | 0 | 0 | $6.8 \cdot 10^8$ | 0 |
| Elastic solid (Copper) | 4.22 | 0 | 0 | $3.42 \cdot 10^{10}$ | $5 \cdot 10^{10}$ |

where

$$p(\rho, s) = -p_\infty - a\rho^2 + \kappa(s)\left(\frac{1}{\rho} - b\right)^{-\gamma} \tag{19}$$

Here $\kappa(s) = \exp(\frac{s}{c_v})$ and $c_v$, $\gamma$, $p_\infty$, $a$, $b$, $\chi$ are positive constants that characterize a given material.

This energy function includes several different physical effects. The first term in the energy function accounts for the phenomenon that an isoentropic compression induces an increase of internal energy ($a$ and $b$ correspond to the van der Waals parameters and have the usual interpretation). The term including the parameter $p_\infty$ accounts for the phenomenon that the energy must increase when the density is reduced. This term is necessary to model fluid or solid material where intermolecular forces are present (see for example [9,12]). The last term in the energy expression models a neohookean elastic solid where the constant $\chi$ is the shear elastic modulus. This term typically accounts for deformations that are finite. For very large deformations, more complete models (Mooney–Rivlin, Ogden) take into account an additional invariant in the energy function in 3D. Without loosing in generality, we stick to a neohookean model as it leads to simpler expressions in the following developments. Moreover, in 2D there are only two independent invariants. Classical models are obtained by particular choices of the coefficients (see Table 1).

Other energy functions can be envisaged to mitigate or accentuate certain aspects of the resulting non-linear relation between stress and deformation, see for example Eq. (49) for the more general Mie–Grüneisen equation of state.

## 3. Numerical scheme

Let $x = (x_1, x_2)$ be the coordinates in the canonical basis of $\mathbb{R}^2$, $u = (u_1, u_2)$ the velocity components, $Y^i_{.j}$ the components of the tensor $[\nabla Y]$ and $\sigma^{ij}$ the components of the stress tensor $\sigma$. Eqs. (16) become

$$
\begin{pmatrix}
\rho \\
\phi_1 \\
\phi_2 \\
Y^1_{,1} \\
Y^2_{,1} \\
Y^1_{,2} \\
Y^2_{,2} \\
\psi
\end{pmatrix}_{,t}
+
\begin{pmatrix}
\phi_1 \\
\frac{(\phi_1)^2}{\rho} - \sigma^{11} \\
\frac{\phi_1 \phi_2}{\rho} - \sigma^{21} \\
\frac{\phi_1 Y^1_{,1} + \phi_2 Y^1_{,2}}{\rho} \\
\frac{\phi_1 Y^2_{,1} + \phi_2 Y^2_{,2}}{\rho} \\
0 \\
0 \\
\frac{\phi_1 \psi - (\sigma^{11}\phi_1 + \sigma^{21}\phi_2)}{\rho}
\end{pmatrix}_{,1}
+
\begin{pmatrix}
\phi_2 \\
\frac{\phi_1 \phi_2}{\rho} - \sigma^{12} \\
\frac{(\phi_2)^2}{\rho} - \sigma^{22} \\
0 \\
0 \\
\frac{\phi_1 Y^1_{,1} + \phi_2 Y^1_{,2}}{\rho} \\
\frac{\phi_1 Y^2_{,1} + \phi_2 Y^2_{,2}}{\rho} \\
\frac{\phi_2 \psi - (\sigma^{12}\phi_1 + \sigma^{22}\phi_2)}{\rho}
\end{pmatrix}_{,2}
=
\begin{pmatrix}
0 \\
0 \\
0 \\
0 \\
0 \\
0 \\
0 \\
0
\end{pmatrix}
\tag{20}
$$

where $\phi_i = \rho u_i$ and $\psi = \rho e$. This system can be written in the compact form

$$\Phi_t + \left(G^1(\Phi)\right)_{,1} + \left(G^2(\Phi)\right)_{,2} = 0 \tag{21}$$

We discretize (21) with a finite volume method on a Cartesian mesh. Let $\Delta x_i$ be the grid spacing in the $x_i$ direction and $\Omega_{i,j}$ the control volume centered at the node $(i\Delta x_1, j\Delta x_2)$. The semi-discretization in space of (21) on $\Omega_{i,j}$ gives

$$(\Phi_{i,j})_t + \frac{G^1_{i+1/2,j} - G^1_{i-1/2,j}}{\Delta x_1} + \frac{G^2_{i,j+1/2} - G^2_{i,j-1/2}}{\Delta x_2} = 0 \tag{22}$$

where $\Phi_{i,j}$ is the value of the conservative variable integrated on $\Omega_{i,j}$ (see Fig. 2).

The fluxes in (22) will be computed by approximate 1D Riemann solvers in the direction orthogonal to the cell sides of the Cartesian mesh. Therefore

$$G^1_{i-1/2,j} \approx \mathcal{F}(\Phi_{i-1,j}; \Phi_{i,j}) \qquad G^1_{i+1/2,j} \approx \mathcal{F}(\Phi_{i,j}; \Phi_{i+1,j}) \tag{23}$$

$$G^2_{i,j-1/2} \approx \mathcal{F}(\Phi_{i,j-1}; \Phi_{i,j}) \qquad G^2_{i,j+1/2} \approx \mathcal{F}(\Phi_{i,j}; \Phi_{i,j+1}) \tag{24}$$

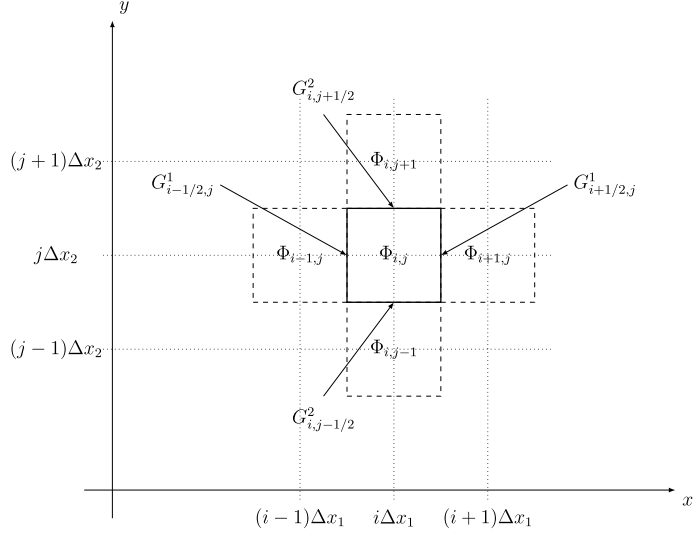

**Fig. 2.** Discretization on the control volume $\Omega_{i,j}$.

where $\mathcal{F}(\cdot,\cdot)$ is a numerical flux function that will be specified in the next sections. The fluxes in (20) have symmetric properties and $(Y_{,2}^1)_t = (Y_{,2}^2)_t = 0$ in the $x_1$ direction. Hence we focus in the next sections on the wave structure of the following one-dimensional problem

$$\Psi_t + \left(F(\Psi)\right)_{,1} = 0 \tag{25}$$

where

$$\Psi = \begin{pmatrix} \rho \\ \phi_1 \\ \phi_2 \\ Y_{,1}^1 \\ Y_{,1}^2 \\ \psi \end{pmatrix} \qquad F(\Psi) = \begin{pmatrix} \phi_1 \\ \dfrac{(\phi_1)^2}{\rho} - \sigma^{11} \\ \dfrac{\phi_1\phi_2}{\rho} - \sigma^{21} \\ \dfrac{\phi_1 Y_{,1}^1 + \phi_2 Y_{,2}^1}{\rho} \\ \dfrac{\phi_1 Y_{,1}^2 + \phi_2 Y_{,2}^2}{\rho} \\ \dfrac{\phi_1\psi - (\sigma^{11}\phi_1 + \sigma^{21}\phi_2)}{\rho} \end{pmatrix}$$

## 4. Characteristic speeds

### 4.1. Equation on entropy

The equations of conservation of momentum and energy (9) give

$$(\rho_0\varepsilon)_t = \left(\rho_0 e - \rho_0\frac{|u|^2}{2}\right)_t = \mathrm{div}_\xi\left(\mathcal{T}^T u\right) - \mathrm{div}_\xi(\mathcal{T})\cdot u = \mathcal{T}:\nabla_\xi u \tag{26}$$

We get from (10)

$$(\rho_0\varepsilon)_t = \rho_0\frac{\partial\varepsilon}{\partial F}:\nabla_\xi u + \rho_0\frac{\partial\varepsilon}{\partial s}\left(s\big(X(\xi,t),t\big)\right)_t$$

If $\frac{\partial\varepsilon}{\partial s} \neq 0$ we get with (11) and (26) the relation $(s(X(\xi,t),t))_t = 0$ or in the Eulerian frame

$$s_t + u\cdot\nabla s = 0 \tag{27}$$

As expected, in this adiabatic and inviscid model entropy is just transported along the characteristics for smooth solutions (no shocks).

## 4.2. Wave speeds in the general case

The governing equations (25) are closed with the constitutive law (13)–(15) which define $\sigma$ as a non-linear function of the unknowns. The wave velocities are locally defined by infinitesimal variation of the conservative variables $\Psi$. Therefore, the energy equation can be simply replaced by (27). Thus we replace (25) by

$$
\begin{pmatrix} \rho \\ \phi_1 \\ \phi_2 \\ Y^1_{,1} \\ Y^2_{,1} \\ s \end{pmatrix}_{,t}
+
\begin{pmatrix}
0 & 1 & 0 & 0 & 0 & 0 \\
-\dfrac{(\phi_1)^2}{\rho^2} & \dfrac{2\phi_1}{\rho} & 0 & -\sigma^{11}_{,1} & -\sigma^{11}_{,2} & -\sigma^{11}_{,s} \\
-\dfrac{\phi_1\phi_2}{\rho^2} & \dfrac{\phi_2}{\rho} & \dfrac{\phi_1}{\rho} & -\sigma^{21}_{,1} & -\sigma^{21}_{,2} & -\sigma^{21}_{,s} \\
-\dfrac{\phi_1 Y^1_{,1}+\phi_2 Y^1_{,2}}{\rho^2} & \dfrac{Y^1_{,1}}{\rho} & \dfrac{Y^1_{,2}}{\rho} & \dfrac{\phi_1}{\rho} & 0 & 0 \\
-\dfrac{\phi_1 Y^2_{,1}+\phi_2 Y^2_{,2}}{\rho^2} & \dfrac{Y^2_{,1}}{\rho} & \dfrac{Y^2_{,2}}{\rho} & 0 & \dfrac{\phi_1}{\rho} & 0 \\
0 & 0 & 0 & 0 & 0 & \dfrac{\phi_1}{\rho}
\end{pmatrix}
\begin{pmatrix} \rho \\ \phi_1 \\ \phi_2 \\ Y^1_{,1} \\ Y^2_{,1} \\ s \end{pmatrix}_{,1}
=
\begin{pmatrix} 0 \\ 0 \\ 0 \\ 0 \\ 0 \\ 0 \end{pmatrix}
\tag{28}
$$

where $\sigma^{ij}_{,1}$, $\sigma^{ij}_{,2}$ and $\sigma^{ij}_{,s}$ denote the derivative of $\sigma^{ij}$ with respect to $Y^1_{,1}$, $Y^2_{,1}$ and $s$ respectively. We denote by $J_{ac}$ the matrix appearing in (28).

We introduce the notations

$$
[\nabla\sigma][\nabla Y] = \begin{pmatrix} \sigma^{11}_{,1} & \sigma^{11}_{,2} \\ \sigma^{21}_{,1} & \sigma^{21}_{,2} \end{pmatrix} \begin{pmatrix} Y^1_{,1} & Y^1_{,2} \\ Y^2_{,1} & Y^2_{,2} \end{pmatrix} = \begin{pmatrix} \sigma^{11}_{,1}Y^1_{,1}+\sigma^{11}_{,2}Y^2_{,1} & \sigma^{11}_{,1}Y^1_{,2}+\sigma^{11}_{,2}Y^2_{,2} \\ \sigma^{21}_{,1}Y^1_{,1}+\sigma^{21}_{,2}Y^2_{,1} & \sigma^{21}_{,1}Y^1_{,2}+\sigma^{21}_{,2}Y^2_{,2} \end{pmatrix}
$$

The characteristic polynomial of $J_{ac}$ can be written under the form

$$
P(\lambda) = \frac{(\lambda - u_1)^2}{\rho^2}\left( \big((\lambda-u_1)^2\rho\big)^2 - \mathrm{Tr}\big(-[\nabla\sigma][\nabla Y]\big)(\lambda-u_1)^2\,\rho + \mathrm{Det}\big(-[\nabla\sigma][\nabla Y]\big) \right)
$$

so that the eigenvalues of $J_{ac}$ are

$$
\Lambda^E = \left\{ u_1, u_1, u_1 + \sqrt{\frac{\alpha_1}{\rho}}, u_1 - \sqrt{\frac{\alpha_1}{\rho}}, u_1 + \sqrt{\frac{\alpha_2}{\rho}}, u_1 - \sqrt{\frac{\alpha_2}{\rho}} \right\}
\tag{29}
$$

where $\alpha_1$ and $\alpha_2$ are the eigenvalues of $-[\nabla\sigma][\nabla Y]$. Therefore, the conditions for the system (25) to be hyperbolic are $\alpha_1 > 0$ and $\alpha_2 > 0$.

## 4.3. Hyperbolicity of the neohookean model

For the neohookean constitutive law (17) we have with (18)

$$
\sigma^{11} = -p(\rho,s) + \chi\big((Y^1_{,2})^2 + (Y^2_{,2})^2 - (Y^1_{,1})^2 - (Y^2_{,1})^2\big) \qquad \sigma^{21} = -2\chi\big(Y^1_{,1}Y^1_{,2} + Y^2_{,1}Y^2_{,2}\big)
$$

where

$$
p(\rho,s) = -p_\infty - a\rho^2 + \kappa(s)\left(\frac{1}{\rho} - b\right)^{-\gamma}
$$

With the relations (29) we find the following explicit expressions for the characteristic speeds

$$
\Lambda^E = \left\{ u_1, u_1, u_1 \pm \sqrt{\frac{\mathcal{A}_1 \pm \sqrt{\mathcal{A}_2}}{\rho}} \right\}
\tag{30}
$$

where

$$
\mathcal{A}_1 = \frac{\rho c^2}{2} + \chi(\alpha_Y + \beta_Y) \qquad \mathcal{A}_2 = \left(\frac{\rho c^2}{2} + \chi(\alpha_Y - \beta_Y)\right)^2 + 4\chi^2(\delta_Y)^2
$$

with

$$
\alpha_Y = (Y^1_{,1})^2 + (Y^2_{,1})^2 \qquad \beta_Y = (Y^1_{,2})^2 + (Y^2_{,2})^2 \qquad \delta_Y = Y^1_{,1}Y^1_{,2} + Y^2_{,1}Y^2_{,2}
$$

and

$$
c^2(\rho,s) = \left.\frac{\partial p}{\partial \rho}\right|_s = \frac{\gamma\kappa(s)}{\rho^2}\left(\frac{1}{\rho} - b\right)^{-\gamma-1} - 2a\rho
\tag{31}
$$

which is the sound speed for a general gas. Hence, when $\chi = 0$ the eigenvalues (30) reduce to the wave speeds of a general gas $\{u_1, u_1, u_1 \pm c, u_1 \pm c\}$. If $\chi > 0$, under the hypothesis that $c^2(\rho,s) > 0$, we have that $\mathcal{A}_1 \geq \sqrt{\mathcal{A}_2} \geq 0$ and the system (25) is hyperbolic for the neohookean model (17).

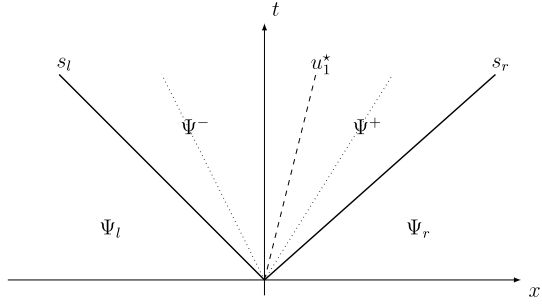

**Fig. 3.** HLLC solver wave pattern.

## 5. Numerical flux

Here we first introduce the HLLC solver employed for a single material and then we specify the numerical flux for a material discontinuity.

### 5.1. HLLC solver

We consider Eq. (25) with the initial condition

$$\Psi(x, t = 0) = \begin{cases} \Psi_l & \text{if } x \leq 0 \\ \Psi_r & \text{if } x > 0. \end{cases} \tag{32}$$

The numerical flux function $\mathcal{F}(\Psi_l; \Psi_r)$ at the cell interface $x = 0$ is determined based on the solution of the HLLC [31] approximate Riemann problem, similarly to what is done in [12]. Even though the exact wave pattern involves five distinct waves (see (29)), the approximate solver approaches the solution using three waves (the contact discontinuity $u_1^\star$, the fastest leftward and rightward waves $s_l$ and $s_r$), thus defining two intermediate states $\Psi^-$ and $\Psi^+$ (see Fig. 3).

We define

$$\Psi^- = \begin{pmatrix} \rho^- \\ \rho^- u_1^\star \\ \rho^- u_2^- \\ (Y_{,1}^1)^- \\ (Y_{,1}^2)^- \\ \psi^- \end{pmatrix} \qquad \mathcal{F}^- = \begin{pmatrix} \rho^- u_1^\star \\ \rho^- (u_1^\star)^2 - (\sigma^{11})^\star \\ \rho^- u_1^\star u_2^- - (\sigma^{21})^\star \\ u_1^\star (Y_{,1}^1)^- + u_2^- (Y_{,2}^1)^- \\ u_1^\star (Y_{,1}^2)^- + u_2^- (Y_{,2}^2)^- \\ u_1^\star \psi^- - ((\sigma^{11})^\star u_1^\star + (\sigma^{21})^\star u_2^-) \end{pmatrix} \tag{33}$$

$$\Psi^+ = \begin{pmatrix} \rho^+ \\ \rho^+ u_1^\star \\ \rho^+ u_2^+ \\ (Y_{,1}^1)^+ \\ (Y_{,1}^1)^+ \\ \psi^+ \end{pmatrix} \qquad \mathcal{F}^+ = \begin{pmatrix} \rho^+ u_1^\star \\ \rho^+ (u_1^\star)^2 - (\sigma^{11})^\star \\ \rho^+ u_1^\star u_2^+ - (\sigma^{21})^\star \\ u_1^\star (Y_{,1}^1)^+ + u_2^+ (Y_{,2}^1)^+ \\ u_1^\star (Y_{,1}^2)^+ + u_2^+ (Y_{,2}^2)^+ \\ u_1^\star \psi^+ - ((\sigma^{11})^\star u_1^\star + (\sigma^{21})^\star u_2^+) \end{pmatrix} \tag{34}$$

Here the stress $(\sigma^{11})^\star$ and $(\sigma^{21})^\star$ are independent unknowns. The numerical flux at the cell interface $x = 0$ is then given by

$$\mathcal{F}(\Psi_l; \Psi_r) = \begin{cases} F(\Psi_l) & \text{if } 0 \leq s_l \\ \mathcal{F}^- & \text{if } s_l \leq 0 \leq u_1^\star \\ \mathcal{F}^+ & \text{if } u_1^\star \leq 0 \leq s_r \\ F(\Psi_r) & \text{if } s_r \leq 0 \end{cases} \tag{35}$$

where $F$ is defined in (25).

The HLLC scheme is based on the assumption that every wave is a shock and therefore Rankine–Hugoniot relations give

$$\begin{cases} F(\Psi_r) - \mathcal{F}^+ = s_r (\Psi_r - \Psi^+) \\ \mathcal{F}^+ - \mathcal{F}^- = u_1^\star (\Psi^+ - \Psi^-) \\ \mathcal{F}^- - F(\Psi_l) = s_l (\Psi^- - \Psi_l) \end{cases} \tag{36}$$

We introduce the notations

$$Q_l = F(\Psi_l) - s_l \Psi_l \qquad Q_r = F(\Psi_r) - s_r \Psi_r \tag{37}$$

and we denote by $Q_l^i$ (respectively $Q_r^i$) the $i$-th component of $Q_l$ (respectively $Q_r$). The Rankine–Hugoniot relations (36) give

$$u_1^\star = \frac{Q_l^2 - Q_r^2}{Q_l^1 - Q_r^1} \qquad (\sigma^{11})^\star = \frac{Q_l^2 Q_r^1 - Q_l^1 Q_r^2}{Q_l^1 - Q_r^1} \tag{38}$$

$$\rho^- = \frac{Q_l^1}{u_1^\star - s_l} \qquad \rho^+ = \frac{Q_r^1}{u_1^\star - s_r} \tag{39}$$

$$(Y_{,1}^1)^- = \frac{Q_l^4 - u_2^\star (Y_{,2}^1)^-}{u_1^\star - s_l} \qquad (Y_{,1}^1)^+ = \frac{Q_r^4 - u_2^\star (Y_{,2}^1)^+}{u_1^\star - s_r} \tag{40}$$

$$(Y_{,1}^2)^- = \frac{Q_l^5 - u_2^\star (Y_{,2}^2)^-}{u_1^\star - s_l} \qquad (Y_{,1}^2)^+ = \frac{Q_r^5 - u_2^\star (Y_{,2}^2)^+}{u_1^\star - s_r} \tag{41}$$

$$\psi^- = \frac{Q_l^6 + (\sigma^{11})^\star u_1^\star + (\sigma^{21})^\star u_2^\star}{u_1^\star - s_l} \qquad \psi^+ = \frac{Q_r^6 + (\sigma^{11})^\star u_1^\star + (\sigma^{21})^\star u_2^\star}{u_1^\star - s_r} \tag{42}$$

In the case where $\Psi_l$ and $\Psi_r$ are solid states ($\chi \neq 0$), both velocity components and the normal stress are continuous across the contact discontinuity. Hence

$$u_2^- = u_2^+ = u_2^\star = \frac{Q_l^3 - Q_r^3}{Q_l^1 - Q_r^1} \qquad (\sigma^{21})^\star = \frac{Q_l^3 Q_r^1 - Q_l^1 Q_r^3}{Q_l^1 - Q_r^1} \tag{43}$$

In the case where at least one of the states $\Psi_l$ or $\Psi_r$ is fluid ($\chi = 0$), $\sigma^{21}$ vanishes at the interface and therefore the transverse velocity can be discontinuous. Hence

$$u_2^- = \frac{Q_l^3}{Q_l^1} \qquad u_2^+ = \frac{Q_r^3}{Q_r^1} \qquad (\sigma^{21})^\star = 0 \tag{44}$$

From Eqs. (36) we consistently obtain $(Y_{,2}^i)^\pm = \frac{(Y_{,2}^i)_l + (Y_{,2}^i)_r}{2}$. This is different compared to [9].

Finally, the robustness of the scheme is strongly influenced by the estimation of $s_l$ and $s_r$. We use the estimate presented in [7] which is a simple way to obtain robust speed estimates:

$$s_l = \min\big((u_1 - \lambda)_l, (u_1 - \lambda)_r\big) \qquad s_r = \max\big((u_1 + \lambda)_l, (u_1 + \lambda)_r\big),$$

where $\lambda = \sqrt{\frac{\mathcal{A}_1 + \sqrt{\mathcal{A}_2}}{\rho}}$, see (30).

### 5.2. Multimaterial solver

The multimaterial solver is detailed in 1D for sake of clarity. In 2D we use exactly the same method in both directions. We consider a case where the material discontinuity is between the cell centers $k-1$ and $k$. The material discontinuity can separate materials with different constitutive laws or discontinuous initial states of the same material.

The main idea of the multimaterial solver is that, instead of (35) we take

$$\mathcal{F}_{k-1/2}^l = \mathcal{F}^- \qquad \mathcal{F}_{k-1/2}^r = \mathcal{F}^+ \tag{45}$$

see Fig. 4 and Eqs. (33) and (34).

As for ghost-fluid methods, the scheme is locally non-conservative since $\mathcal{F}^- \neq \mathcal{F}^+$, but it is consistent since $\mathcal{F}^\pm$ are regular enough functions of the states to the left and to the right of the interface and $\mathcal{F}^+ = \mathcal{F}^-$ when those states are identical. As shown in the numerical test section the error in conservation is negligible: the shock speeds and positions are correctly predicted. Indeed, the number of cell interfaces for which a non-conservative numerical flux is employed is always negligible compared to the total number of mesh cells.

A special care is needed at the multimaterial interface. To avoid the singularity due to the discontinuity of $Y$ and to keep the interface sharp, the computation of $(Y_{,2}^i)^\pm$ is one-sided:

$$(Y_{,2}^i)^- = (Y_{,2}^i)_l \qquad (Y_{,2}^i)^+ = (Y_{,2}^i)_r \tag{46}$$

This can be seen as an extension by continuity of $Y_{,2}^i$ relative to each material.

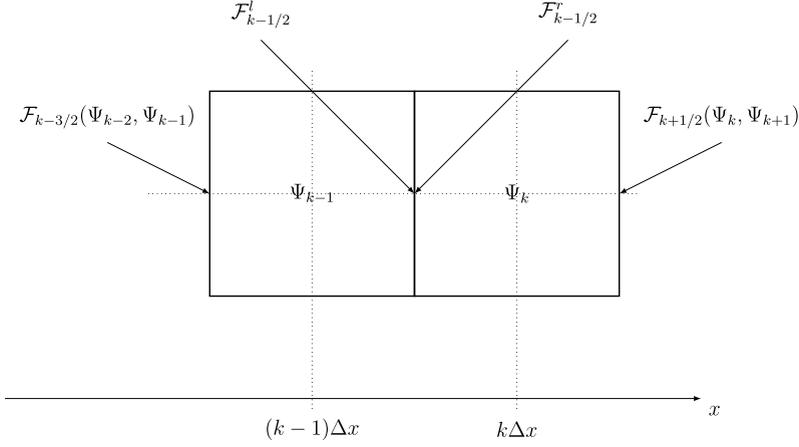

**Fig. 4.** Fluxes at the material discontinuity.

*5.3. Extension to second-order in each space direction*

The scheme is extended to second-order accuracy in space with a piecewise-linear slope reconstruction (MUSCL). For the cells separated by the multimaterial interface, the stencil used to compute the slopes is smaller. We calculate these slopes using the corresponding intermediate state of the multimaterial Riemann problem. For example, if cells $k-1$ and $k$ do not belong to the same medium, the slope of cell $k-1$ is computed using $\Psi_{k-2}$, $\Psi_{k-1}$ and the intermediate state to the left of the contact discontinuity. This state is determined thanks to the solution of an approximate Riemann problem between $\Psi_{k-1}$ and $\Psi_k$, without slope reconstruction.

## 6. Interface advection and time integration

Coherently with the fully Eulerian approach, a level set function is used to follow the interface separating the two materials. The level-set function is transported with the velocity field by the equation:

$$\varphi_t + u \cdot \nabla \varphi = 0 \tag{47}$$

This equation is approximated with a WENO 5 scheme [18].

The conservation equations and the interface advection are explicitly integrated in time by a Runge–Kutta 2 scheme. The interface position is advected using the material velocity field. For numerical stability, the integration step is limited by the fastest characteristics over the grid points. Hence, the interface position will belong to the same interval between two grid points for more than one time step. When the physical interface overcomes a grid point, the corresponding conservative variables, say $\Psi_k$, do not correspond anymore to the material present at that grid point before the integration step. When the physical interface moves to the right, then we take $\Psi_k = \Psi^-$, whereas if it moves to the left $\Psi_k = \Psi^+$.

Compared to the ghost-fluid method and its improvements and variants [2,8,10], this scheme is simpler as it does not require the storage of any additional variables or equation of state relative to a ghost fluid to treat the material interface, nor the solution of additional Riemann problems, each relative to a different material at the interface. Our scheme only relies on the intermediate states $\Psi^-$ and $\Psi^+$ that are anyway computed at every cell interface.

## 7. Exact Riemann solver for 1D test cases

A five-wave exact Riemann solver based on the method described in [5] is set up in order to validate the numerical scheme in one space dimension. This is an iterative method whose solution converges to the solution of the Riemann problem, i.e. system of equation $\Psi_t + (F(\Psi))_{,1} = 0$ with the initial conditions

$$\Psi(x, t = 0) = \begin{cases} \Psi_l & \text{if } x \leq x_0 \\ \Psi_r & \text{if } x > x_0, \end{cases} \tag{48}$$

where $x_0$ is the position of the discontinuity at initial time. The solution of the Riemann problem is composed of six constant states separated by five distinct waves, which are from left to right: a normal stress wave, a transverse shear wave, a contact wave, a transverse shear wave and a normal stress wave, see Fig. 5. The solution across each wave is uniquely determined knowing the state on either side, the wave type and the wave speed.

Given an initial estimate of each constant state, the wave types and the wave speeds are determined. Then, knowing the wave speeds, the left and right states $\Psi_l$ and $\Psi_r$, and the wave types, one can determine the state across waves 1 and 4 using a Newton method, and so on for waves 2 and 3. Once a solution is found for states $\Psi_2$ and $\Psi_3$, the error concerning the continuity of the velocity and the stress tensor across the contact wave is computed. This error is used to improve wave speeds estimates perturbing each wave speed, see [5] for a full description.

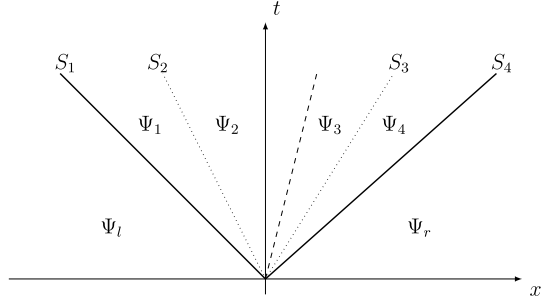

**Fig. 5.** Representation of the Riemann problem for system (16). The constant states are denoted by $\Psi_j$, and the wave speeds concerning the non-linear waves by $S_j$, $1 \le j \le 4$. The dashed line is the contact discontinuity wave.

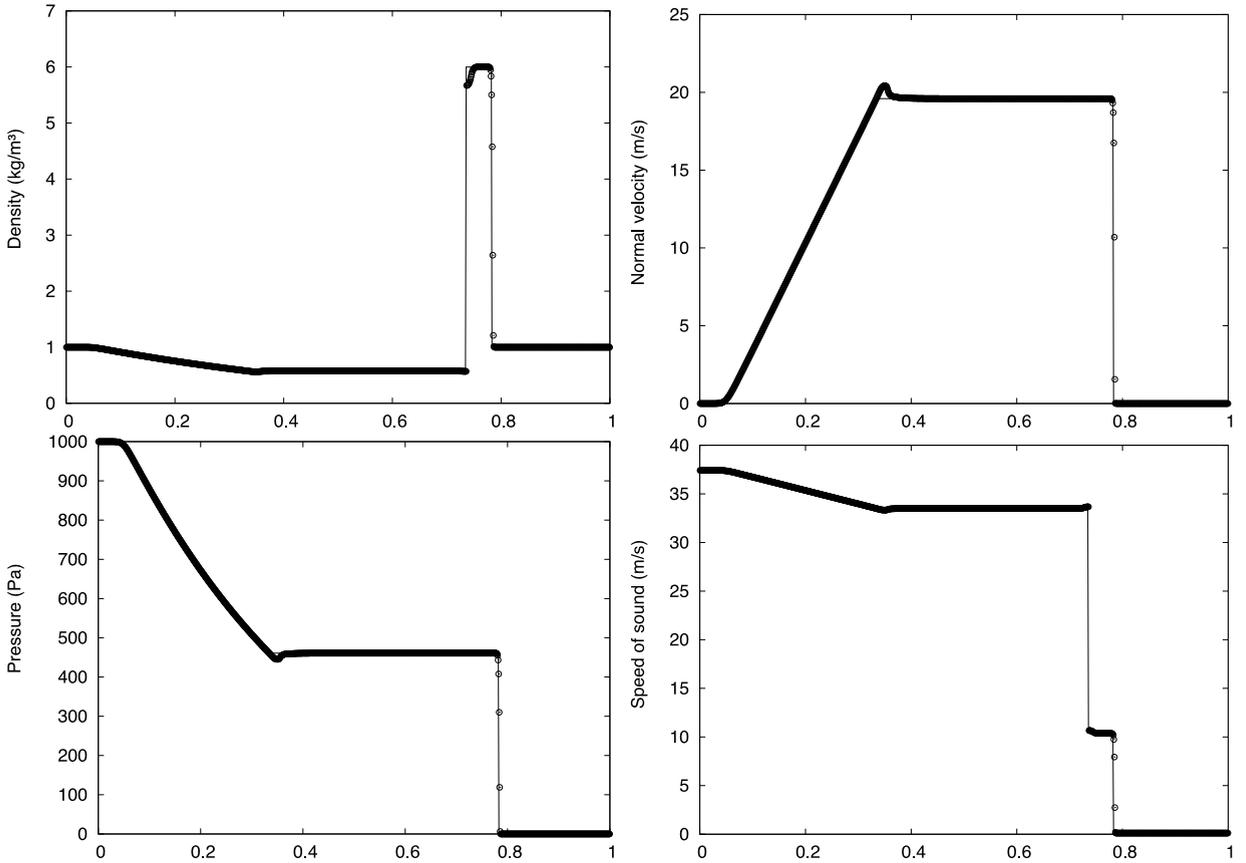

**Fig. 6.** Gas–gas shock tube with high pressure ratio (TC1).

## 8. Numerical results

All the computations are performed with the second order scheme as described above and a minmod slope limiter. First, we present one-dimensional validations with respect to the exact solution and then, two-dimensional shock bubble interactions and impacts.

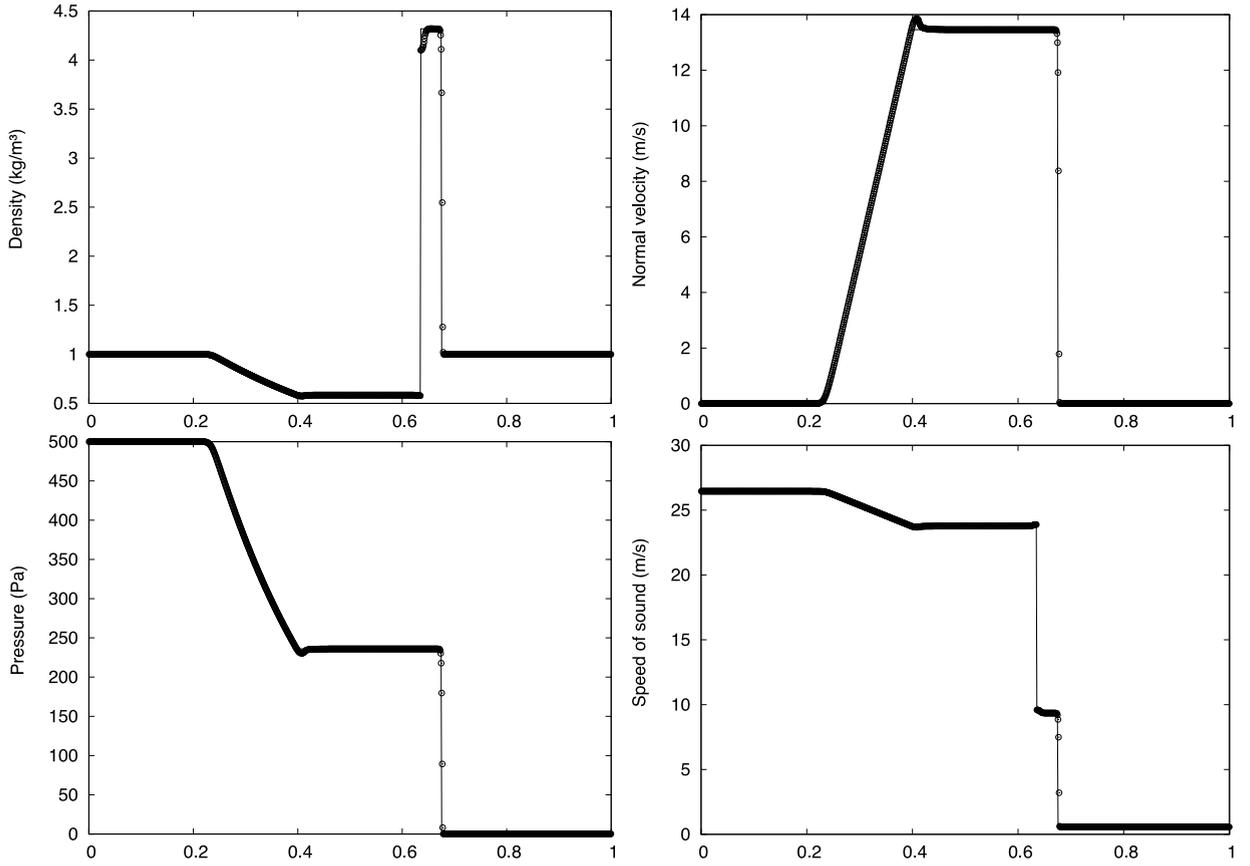

**Fig. 7.** Gas–gas shock tube with high pressure ratio and different constitutive laws (TC2).

**Table 2**
One-dimensional test case description.

| TC | side | $\rho$ [kg/m$^3$] | $u_1$ [m/s] | $u_2$ [m/s] | $p$ [Pa] | $\gamma$ | $p_\infty$ [Pa] | $\chi$ [Pa] | $x_0$ [m] | $t_{end}$ [s] |
|---|---|---|---|---|---|---|---|---|---|---|
| 1 | left | 1 | 0 | 0 | 1000 | 1.4 | 0 | 0 | 0.5 | 0.012 |
|   | right | 1 | 0 | 0 | 0.01 | 1.4 | 0 | 0 | | |
| 2 | left | 1 | 0 | 0 | 500 | 1.4 | 0 | 0 | 0.5 | 0.01 |
|   | right | 1 | 0 | 0 | 0.2 | 1.6 | 0 | 0 | | |
| 3 | left | 1000 | 0 | 0 | $10^9$ | 4.4 | $6.8 \cdot 10^8$ | 0 | 0.7 | $2.4 \cdot 10^{-4}$ |
|   | right | 50 | 0 | 0 | $10^5$ | 1.4 | 0 | 0 | | |
| 4 | left | 8900 | 0 | 0 | $10^9$ | 4.22 | $3.42 \cdot 10^{10}$ | $5 \cdot 10^{10}$ | 0.5 | $5 \cdot 10^{-5}$ |
|   | right | 8900 | 0 | 100 | $10^5$ | 4.22 | $3.42 \cdot 10^{10}$ | $5 \cdot 10^{10}$ | | |
| 5 | left | 8900 | 1000 | 100 | $10^5$ | 4.22 | $3.42 \cdot 10^{10}$ | $5 \cdot 10^{10}$ | 0.5 | $1.5 \cdot 10^{-4}$ |
|   | right | 1 | 1000 | 0 | $10^5$ | 1.4 | 0 | 0 | | |
| 6 | left | 8900 | 0 | 0 | $5 \cdot 10^9$ | 4.22 | $3.42 \cdot 10^{10}$ | $5 \cdot 10^{10}$ | 0.6 | $8.7 \cdot 10^{-5}$ |
|   | right | 50 | 0 | 0 | $10^5$ | 1.4 | 0 | 0 | | |

## 8.1. One-dimensional test cases

We present fluid–fluid, solid–solid and solid–fluid shock tubes problems. The computational domain is always $[0, 1]$, the discontinuity is initially located at $x_0$ and $t_{end}$ is the final time. The initial condition and the physical parameters of the test cases are described in Table 2. In all the figures the circles represent the numerical solution, the solid line the exact solution. The computations are performed on 1000 grid points unless otherwise stated and CFL is 0.6.

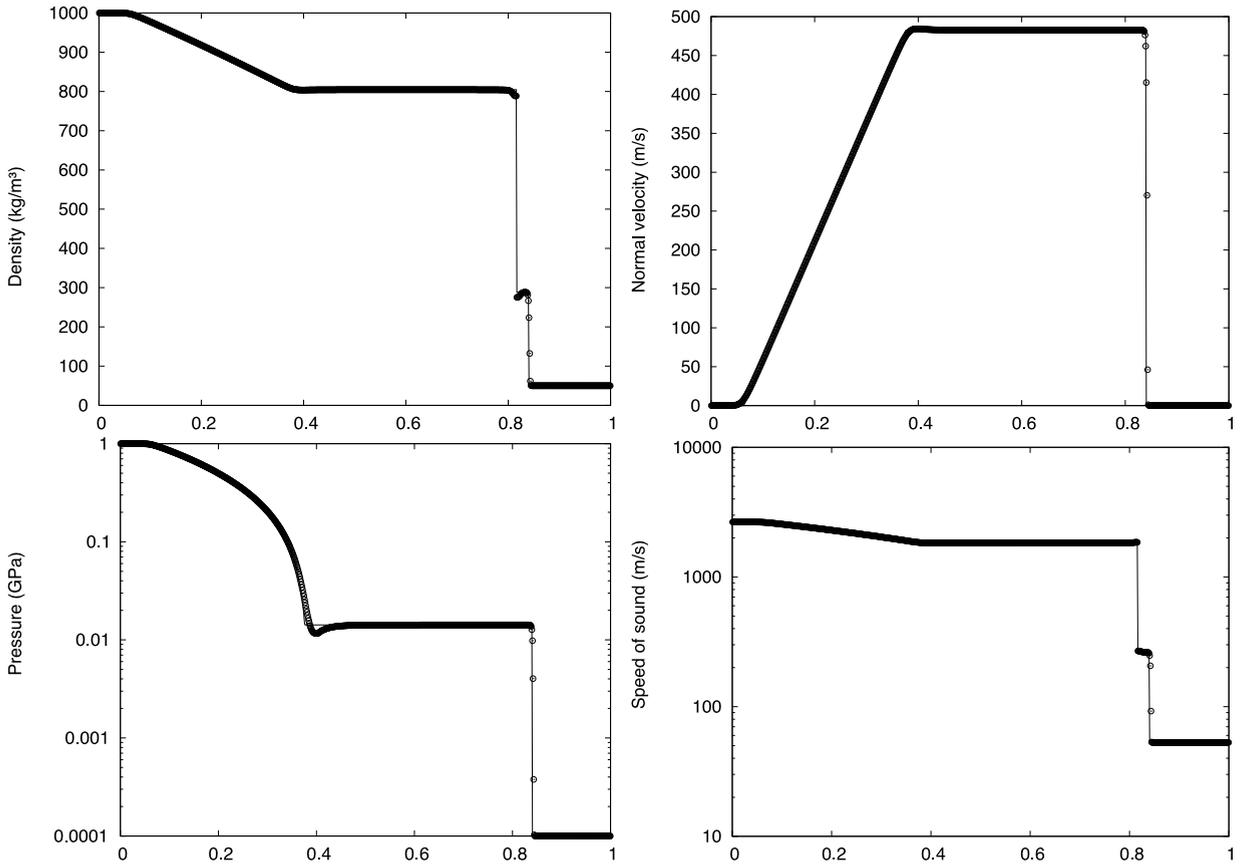

**Fig. 8.** Water–air shock tube with large density and pressure ratios (TC3). On the pressure plot, the scale is logarithmic.

### 8.1.1. Fluid–fluid

The test cases 1 and 2 are shock tubes filled with perfect gases with large pressure ratios, and for test case 2, different adiabatic constants $\gamma$. The results are shown on Figs. 6 and 7. The interface stays sharp, the results are oscillation-free, in good agreement with the exact solution, even if they present the classical overheating problem close to the interface.

The test case 3 is a water–air shock tube. The constitutive laws are different, the density and pressure ratios are large. The results are shown on Fig. 8. Again, the results are in good agreement with the exact solution, are oscillation free, have the correct shock strength and speed, and the interface is kept sharp.

For tests cases 2 and 3, we show density plots for lower grid resolutions: 200 and 500 grid points, see Fig. 9. For these cases, density is typically the most critical variable in terms of discrepancy with respect to the exact solution. The solution is less accurate compared to the finer grid, but the contact discontinuity is sharp and the shocks are smeared but at the right place. In this sense, for test case 3 we present the mass conservation error as a function of time for a coarse grid of 200 points, see Fig. 10. For the same case, we show the convergence of the mass conservation error in the max norm as a function of the grid refinement. It can be observed that mass conservation is consistently granted with a convergence of order 1.

### 8.1.2. Solid–solid

The test case 4 is a solid shock tube with shear. An interface separates the high pressure chamber on the left where the copper is at rest and the same material on the right at low pressure. A tangential velocity discontinuity is imposed, so five waves are expected to appear. The results are shown on Fig. 11.

One can distinguish five waves in the field. The fastest waves are those relative to the normal stress, the middle one is the interface, and the two intermediate waves are those relative to the tangential (shear) stress. The shear and the normal stress, the horizontal and the vertical speeds are continuous at the interface, whereas the density and the pressure are discontinuous at the interface. The fields are oscillation free, in good agreement with the exact solution, and the interface is sharp.

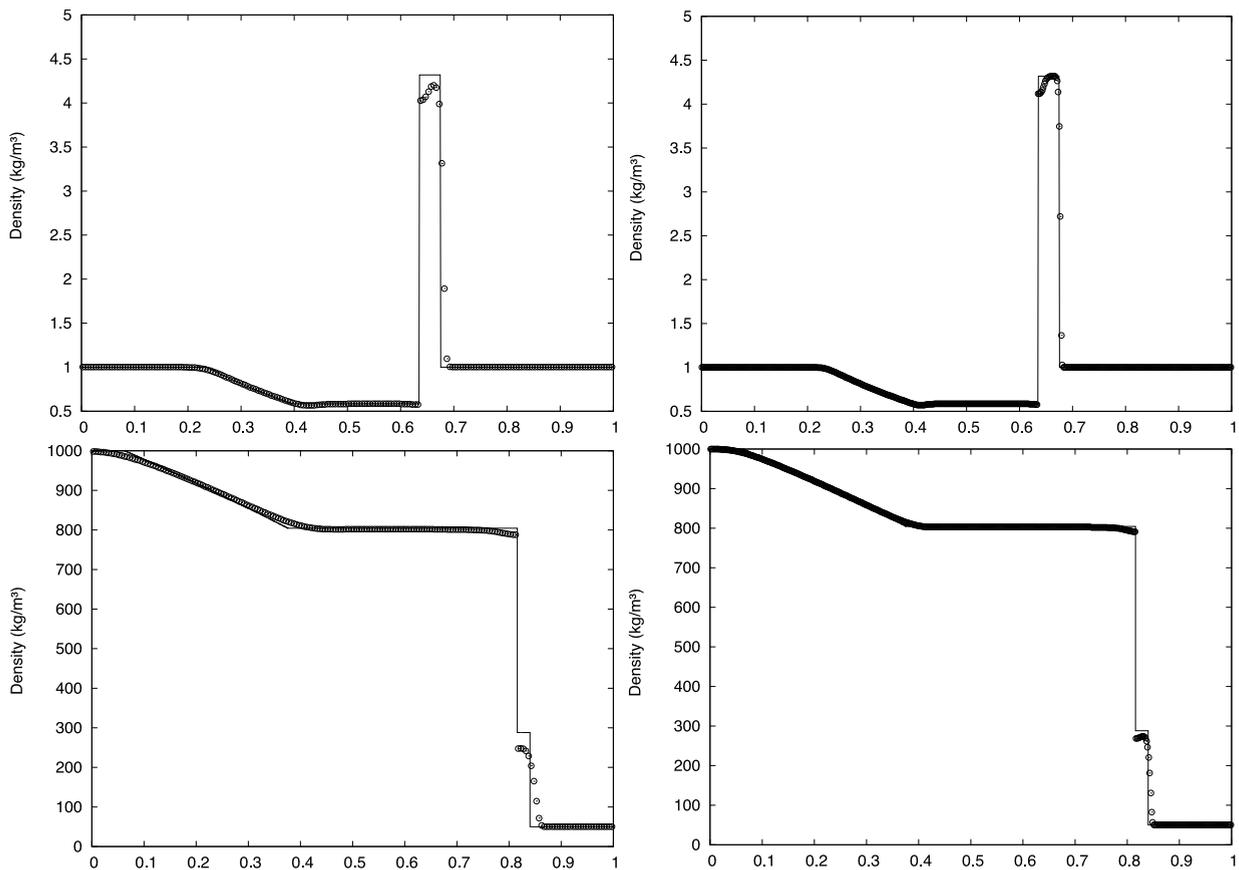

**Fig. 9.** Density plots with 200 (left) and 500 grid points (right). Upper row TC2, lower TC3.

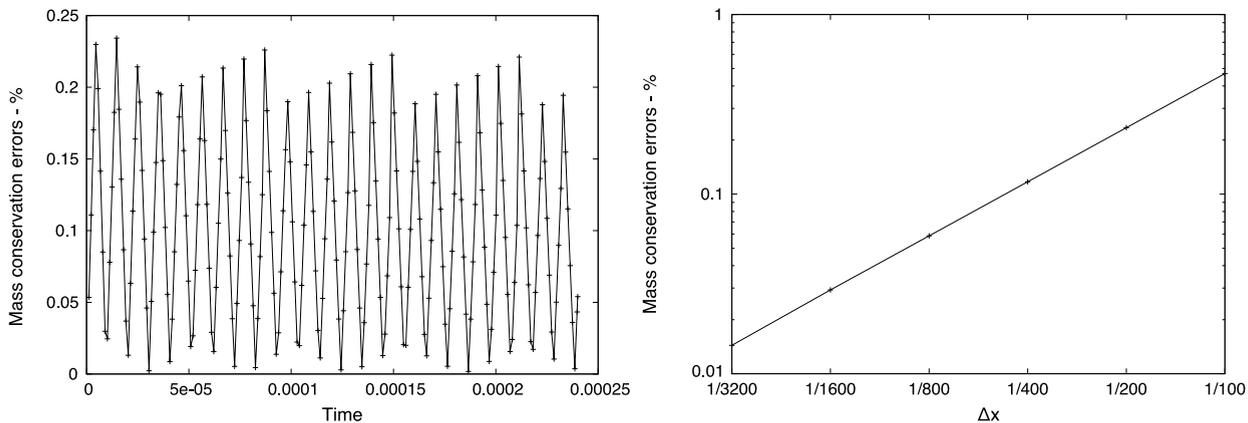

**Fig. 10.** Mass conservation error in percentage of the total mass as a function of time for 200 grid points (left). Convergence of the mass conservation error in the max norm as a function of the grid refinement (right).

### 8.1.3. Solid–fluid

Test case 5 is the advection of a copper–air interface at uniform speed with a tangential velocity discontinuity. This computation is performed on 100 points only because the pattern of the solution is very simple. The results for TC5 are presented on Fig. 12. One can see that the interface is perfectly transported. The flow is oscillation free, and the mechanical equilibrium is preserved in presence of a large density ratio and a transverse velocity discontinuity.

Test case 6 is a shock tube containing copper at high pressure and air at atmospheric pressure. This test case is stiff because at the initial time, the pressure and density ratios are very large. One can see the results on Fig. 13. The normal velocity and the normal stress are continuous at the interface. The pressure and the density are discontinuous at the interface which is kept sharp. We can see a wave transmitted in the air on the density and pressure plots. The results confirm the robustness of the scheme.

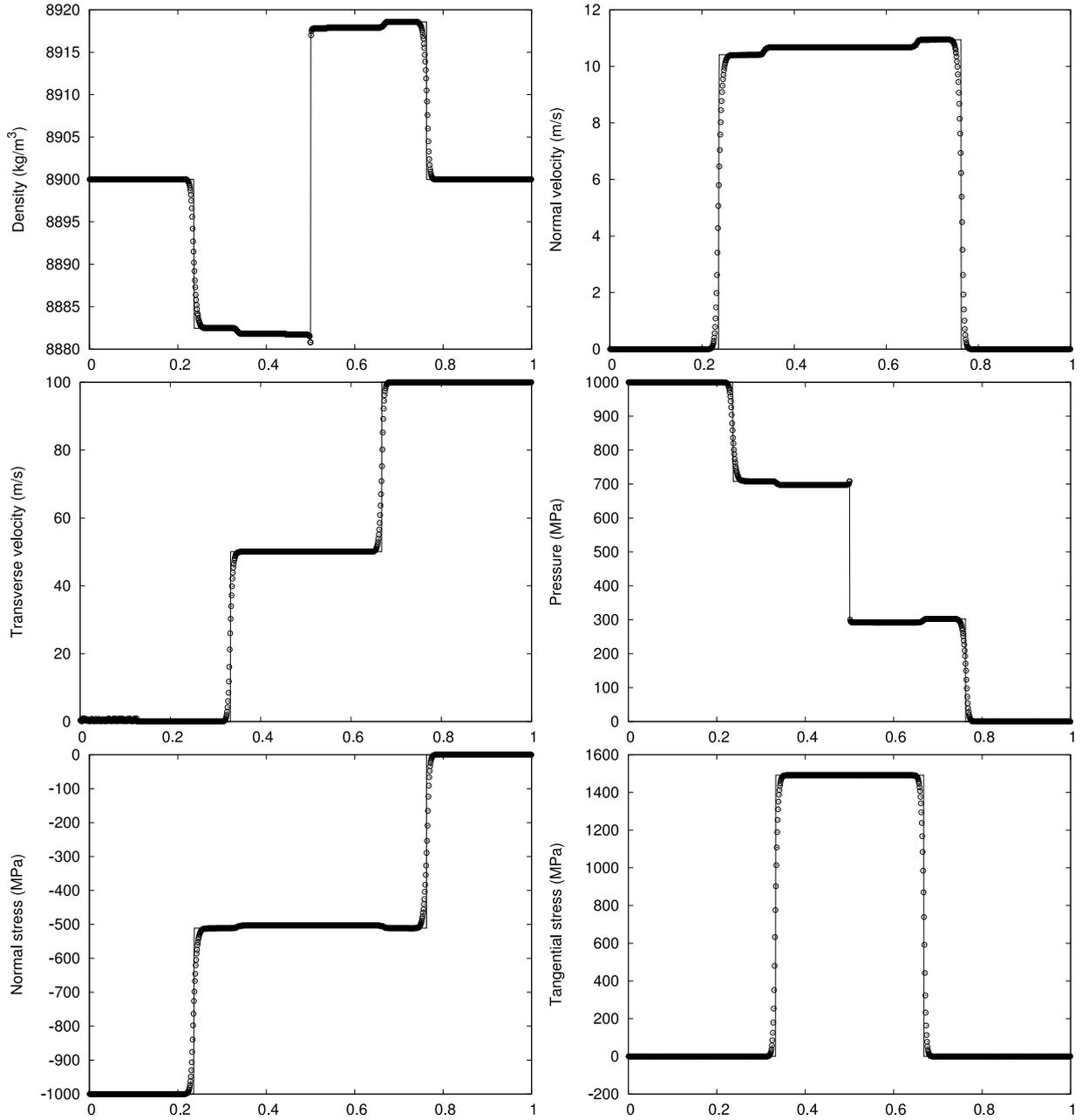

**Fig. 11.** Copper shock tube with shear (TC4).

### 8.1.4. Fluid–fluid with Mie–Grüneisen equation of state

The test case 7 is a shock tube filled with a fluid modeled by a Mie–Grüneisen equation of state. The volumetric part of the energy is written in the form

$$\varepsilon(\rho, s) = \frac{\kappa(s)\rho^{\gamma-1}}{\gamma - 1} + \frac{A_1}{\rho_{ref}(E_1 - 1)}\left(\frac{\rho}{\rho_{ref}}\right)^{E_1-1} - \frac{A_2}{\rho_{ref}(E_2 - 1)}\left(\frac{\rho}{\rho_{ref}}\right)^{E_2-1} \tag{49}$$

The physical parameters are given in Table 3 and the initial conditions in Table 4 are those in [28]. Numerical solutions for density, pressure, velocity are shown in Fig. 14 and compared to the exact solution. The results show that the numerical scheme presented here equally applies to this model.

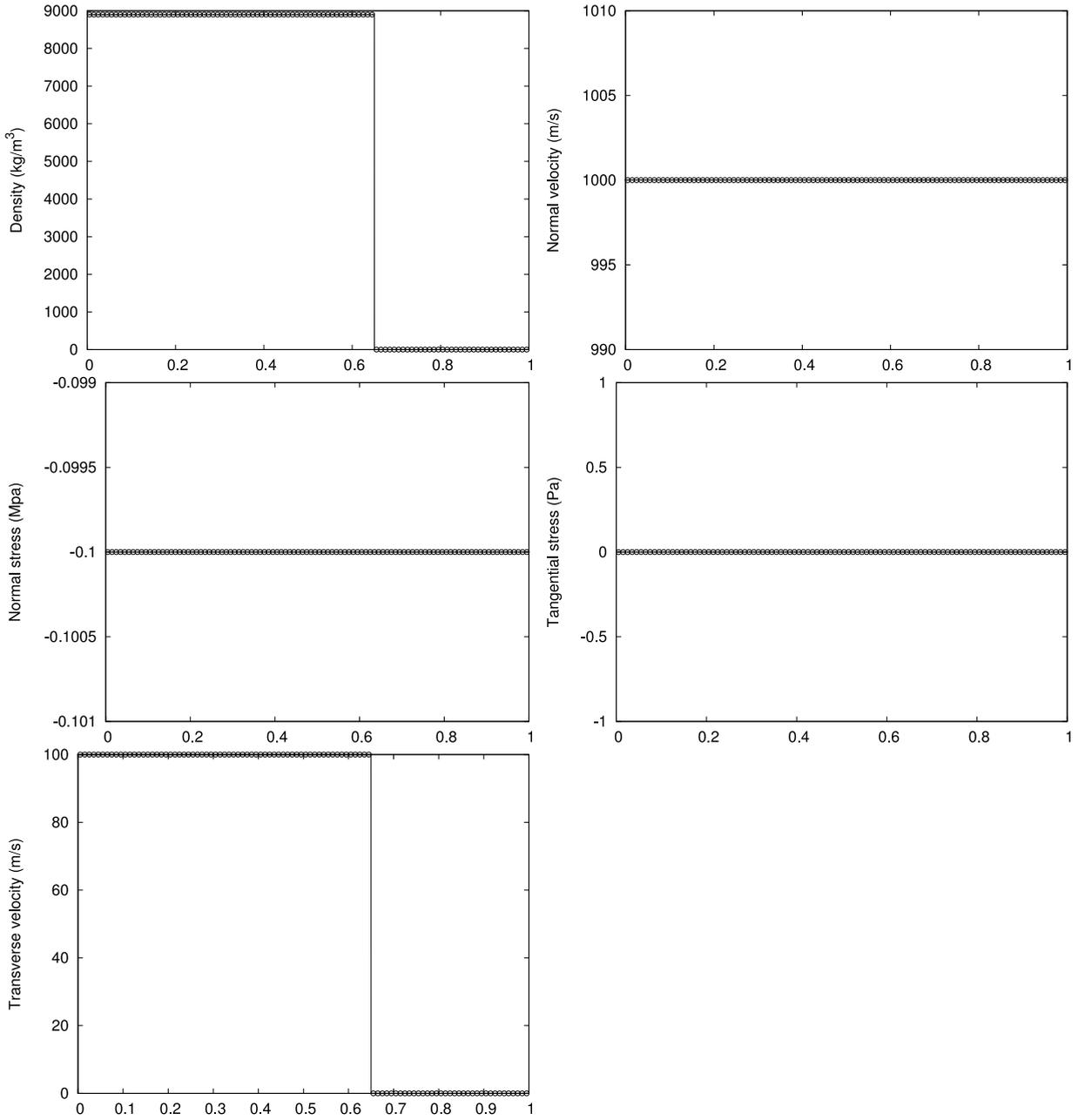

**Fig. 12.** Advection of a copper–air interface with a transverse velocity discontinuity (TC5).

**Table 3**
Parameters for the Mie–Grüneisen equation of state.

| TC | $\gamma$ | $\rho_{ref}$ [kg/m$^3$] | $A_1$ [Pa] | $A_2$ [Pa] | $E_1$ | $E_2$ |
|---|---|---|---|---|---|---|
| 7 | 2.19 | 1134 | $0.819181 \cdot 10^9$ | $1.50835 \cdot 10^9$ | 4.52969 | 1.42144 |

### 8.2. Two-dimensional simulations

Numerical results are presented for the 2D extension of the scheme. Again, computations are performed with the second-order scheme and minmod limiter.

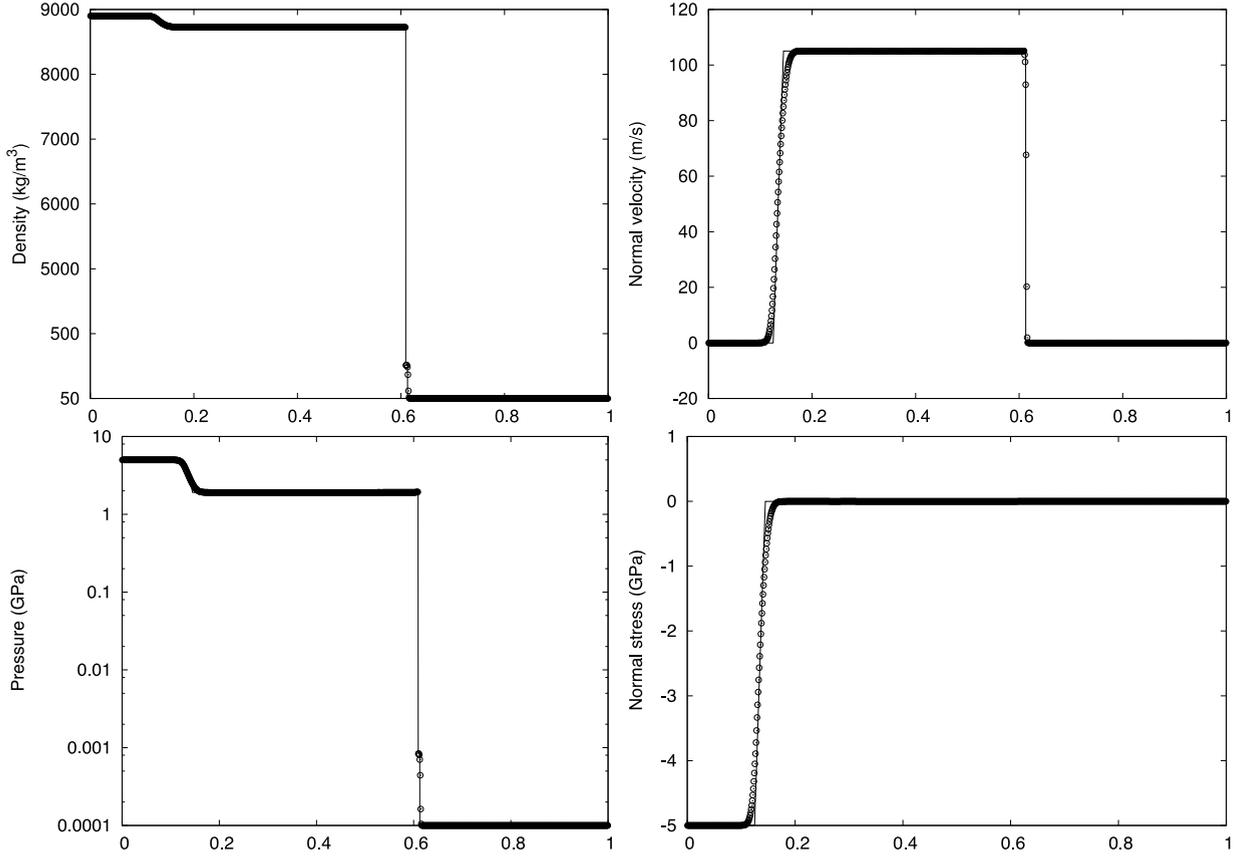

**Fig. 13.** Copper–air shock tube with high pressure and density ratios (TC6). On the density and pressure plot, the scale is logarithmic for the $y$ axis (up to 5000 for the density).

**Table 4**
Mie–Grüneisen test case description.

| TC | side | $\rho$ [kg/m$^3$] | $u_1$ [m/s] | $u_2$ [m/s] | $p$ [Pa] | $x_0$ [m] | $t_{\text{end}}$ [s] |
|---|---|---|---|---|---|---|---|
| 7 | left | 1134 | 0 | 0 | $20 \cdot 10^9$ | 0.6 | $50 \cdot 10^{-6}$ |
| | right | 1200 | 0 | 0 | $0.2 \cdot 10^6$ | | |

### 8.2.1. Shock–bubble interaction

We consider two shock–bubble interaction test cases involving two different fluids. The boundary conditions are reflection on upper and lower borders and homogeneous Neumann conditions for inlet and outlet. Surface tension is not taken into account in the model, so the solution is not physically relevant when the interface presents a high curvature. Hence, the bubble splitting happens when a structure becomes smaller than the size of a computational cell. The computations are performed on a single processor and last for about 60 h of CPU for the finer grids.

*8.2.1.1. Air–helium shock–bubble interaction* The test case 8 is the propagation of a Mach 1.22 shock moving in air, through a helium bubble. This test case has been initially proposed in [26]. The initial configuration and the physical parameters are described in Table 5 and on Fig. 15. The helium bubble is contaminated with 28% of air and we have chosen the initial conditions reported in [22] for comparison purposes.

The computation is performed on a 2000 × 400 points grid. Numerical Schlieren pictures are presented on Fig. 16. These pictures are obtained plotting the $|\nabla \rho|$ field using a logarithmic scale. The evolution of the helium bubble shape is presented on Fig. 17 on a single picture.

One can observe that the shock propagates faster in the helium bubble than in air and that the reflection of the shock inside the bubble presents the typical patterns. Then, the bubble is deformed and the two chambers are linked by a filament which finally breaks. When a structure smaller then a computational cell passes between to two cell centers, the isoline zero of the level set function disappears, and with it, the fluid contained. The two vertical waves in addition to the shock result from the fact that the Rankine–Hugoniot conditions are not exactly satisfied through the shock. From the initial discontinuity arise three waves. The logarithmic scale highlight these waves, but their amplitude is small.

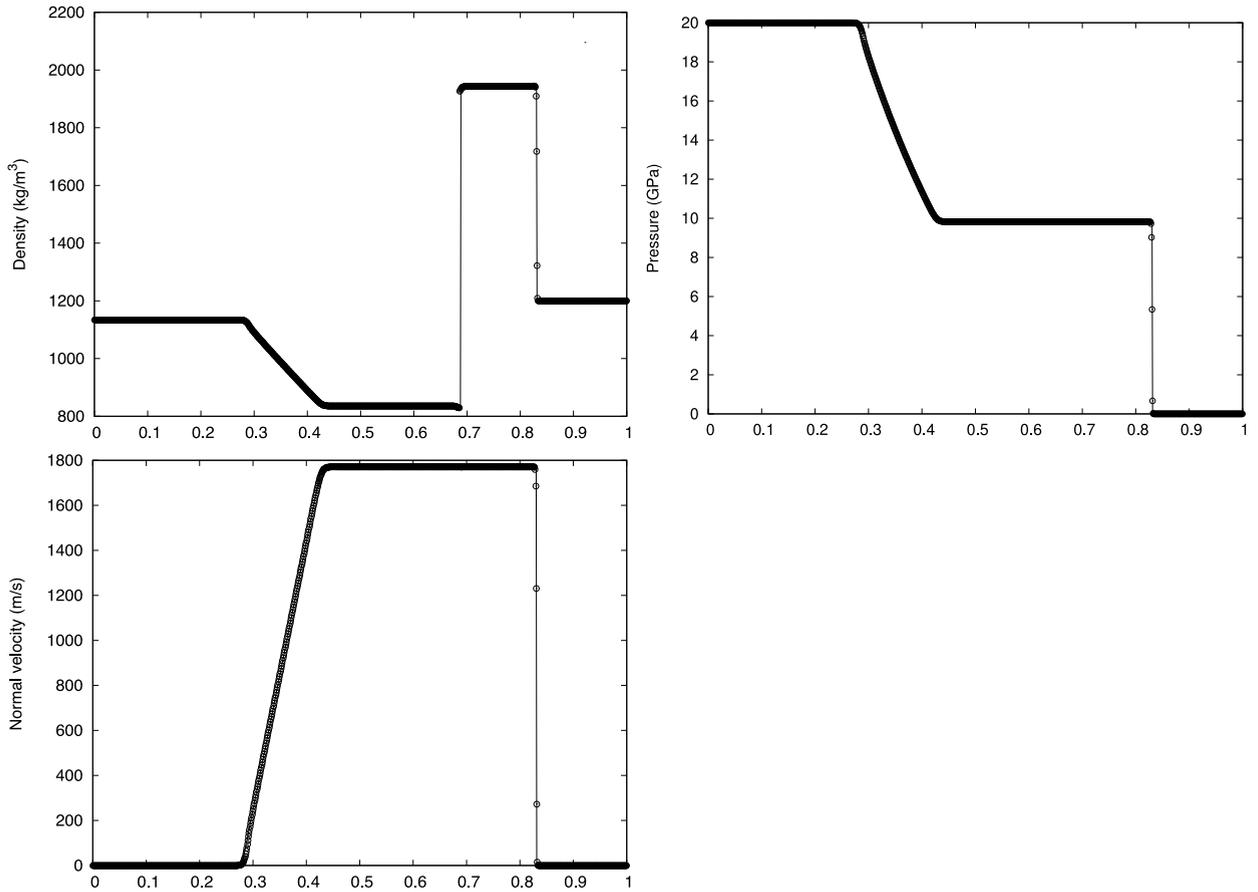

**Fig. 14.** Fluid–fluid shock tube test case with Mie–Grüneisen equation of state (TC7).

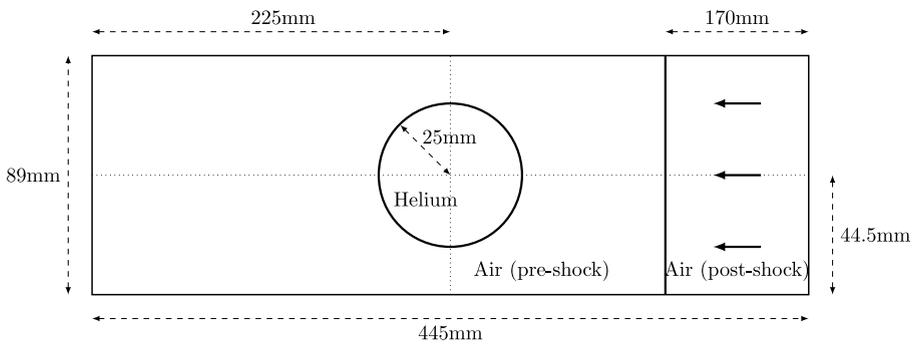

**Fig. 15.** Sketch of the initial configuration for TC8.

This simulation is in overall good qualitative agreement with previous results, see [10,22,26], even though the initial conditions and fluid parameters are not homogeneous through the literature. In particular, there is a remarkable accordance with pictures shown in [22] for corresponding time snapshots, even if our simulations are less resolved.

*8.2.1.2. Air–water shock–bubble interaction* The test case 9 involves a Van der Waals bubble and a Mach 1.422 shock in a stiffened gas. The computational domain is $[-0.2, 1] \times [0, 1]$, the initial configuration and the physical parameters are described in Table 5 and in Fig. 18. This test case is more severe compared to the previous one since it presents larger density ratios.

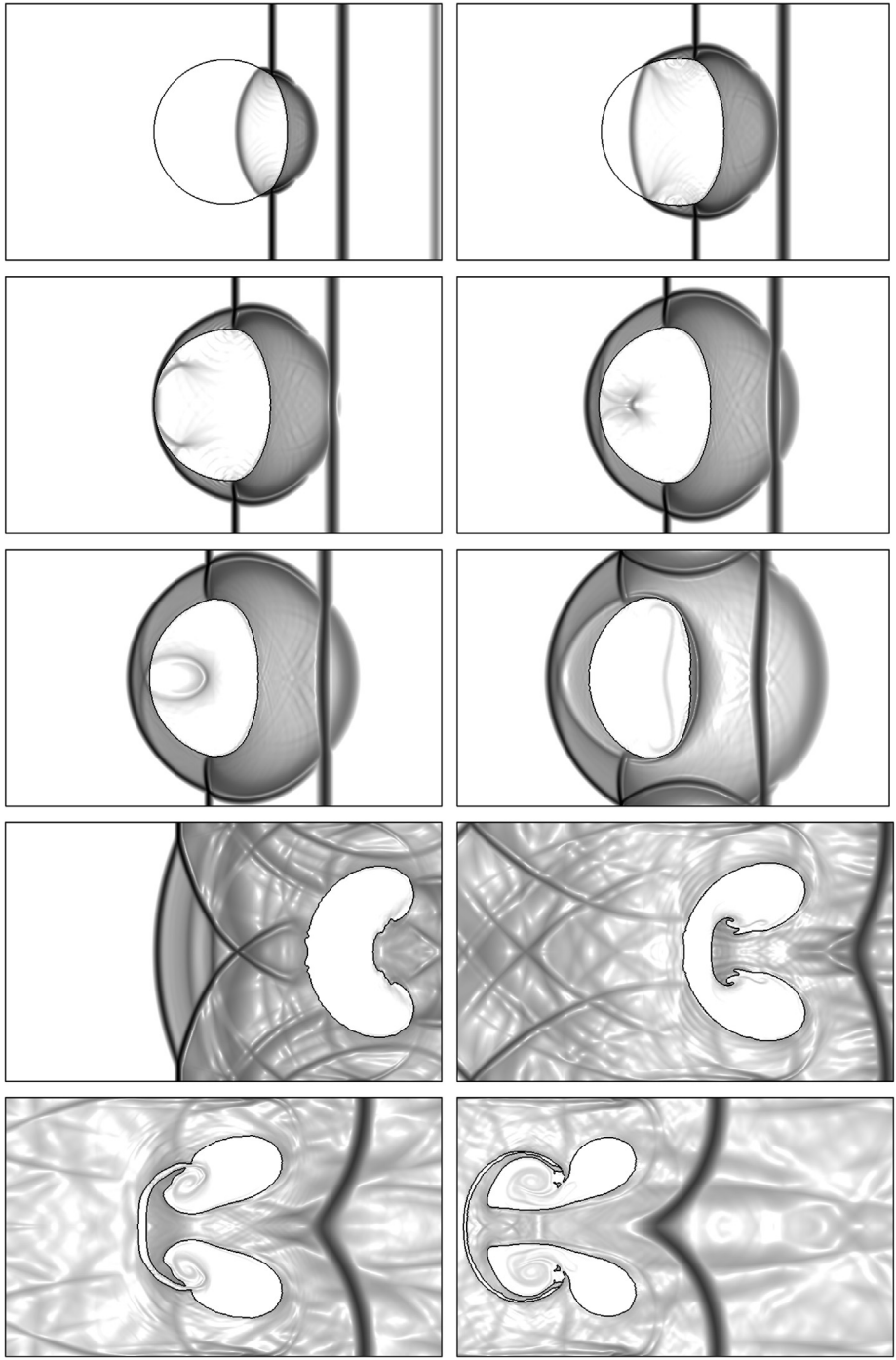

**Fig. 16.** Interaction of a Mach 1.22 shock in air and a helium bubble partially contaminated with 28% of air (TC8). Numerical Schlieren pictures at $t = 23$ μs, 42 μs, 53 μs, 66 μs, 75 μs, 102 μs, 260 μs, 445 μs, 674 μs, 983 μs. From left to right, top to bottom.

The computational grid is $480 \times 400$, and numerical Schlieren pictures as well as pressure and density along the straight line $y = 0.5$ are presented in Figs. 19 and 20.

This simulation is in good agreement with the literature, see [3] and [30] for example. The bubble is strongly compressed, it breaks and swirls. Also, the results compare well with the numerical study in [15] even though a different constitutive law for water is employed in that paper.

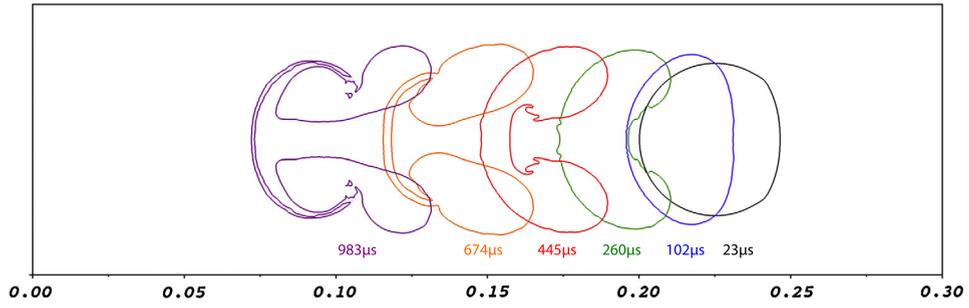

**Fig. 17.** Evolution of the helium bubble shape (TC8).

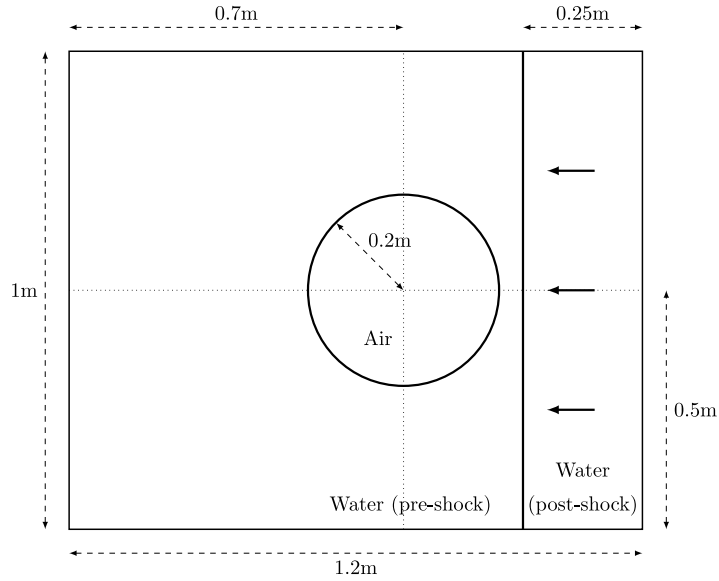

**Fig. 18.** Sketch of the initial configuration for TC9.

**Table 5**
Two-dimensional test case description.

| TC | Media | $\rho$ [kg/m³] | $u_1$ [m/s] | $p$ [Pa] | $\gamma$ | $a$ | $b$ | $p_\infty$ [Pa] | $\chi$ [Pa] |
|----|-------|----------|---------|--------|--------|-----|-----|-----------|---------|
| 8 | Air (pre-shock) | 1.225 | 0 | 101 325 | 1.4 | 0 | 0 | 0 | 0 |
|   | Air (post-shock) | 1.6861 | −113.534 | 159 059 |   |   |   |   |   |
|   | Helium | 0.2228 | 0 | 101 325 | 1.648 | 0 | 0 | 0 | 0 |
| 9 | Water (pre-shock) | 1000 | 0 | $10^5$ | 4.4 | 0 | 0 | $6 \cdot 10^8$ | 0 |
|   | Water (post-shock) | 1230 | −432.69 | $10^9$ |   |   |   |   |   |
|   | Air | 1.2 | 0 | $10^5$ | 1.4 | 5 | $10^{-3}$ | 0 | 0 |
| 10 | Copper (plate) | 8900 | 0 | $10^5$ | 4.22 | 0 | 0 | $3.42 \cdot 10^{10}$ | $5 \cdot 10^{10}$ |
|   | Copper (projectile) | 8900 | 800 | $10^5$ |   |   |   |   |   |
|   | Air | 1 | 0 | $10^5$ | 1.4 | 0 | 0 | 0 | 0 |
| 11 | Liquid (plate) | 8900 | 0 | $10^5$ | 4.22 | 0 | 0 | $3.42 \cdot 10^{10}$ | 0 |
|   | Liquid (projectile) | 8900 | 800 | $10^5$ |   |   |   |   |   |
|   | Air | 1 | 0 | $10^5$ | 1.4 | 0 | 0 | 0 | 0 |

### 8.2.2. Impacts

Finally, we present two impacts simulations of a 800 m/s projectile on a plate in air (TC10 and TC11). This test case is performed in [9] and [20]. The computational domain is $[-0.5, 0.5] \times [-0.5, 0.5]$, the initial configuration and the physical parameters are described in Table 5 and in Fig. 21. The projectile and the plate are adjacent at initial time. Homogeneous Neumann conditions are imposed at the borders. Two computations are performed on a $1000 \times 1000$ mesh, changing the shear parameter ($\chi = 5 \cdot 10^{10}$ in TC10 and $\chi = 0$ in TC11, see Table 5). The computations are performed on a single processor and last for about 50 h of CPU for the finer grids.

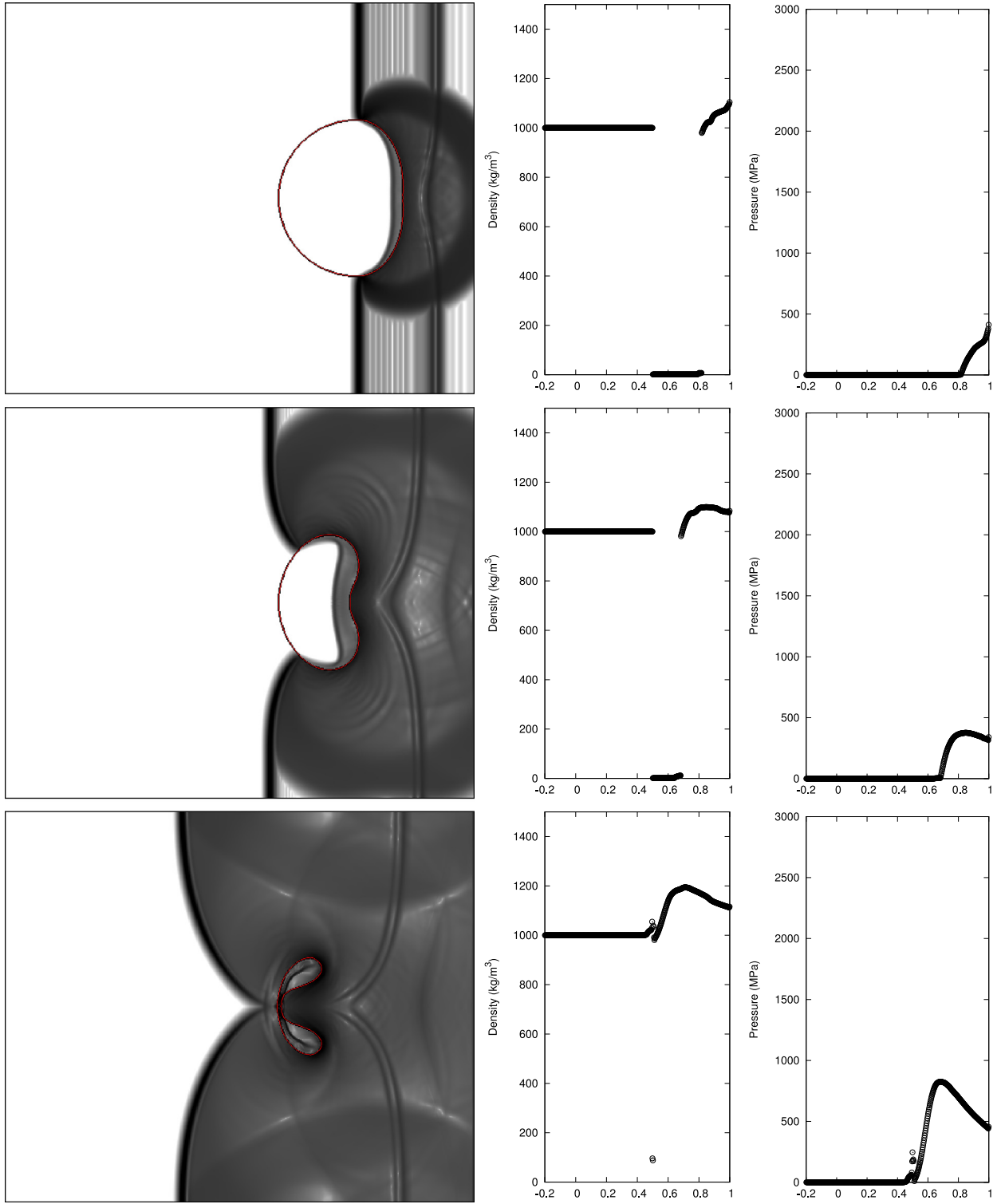

**Fig. 19.** Interaction of a Mach 1.422 water shock and an air bubble (TC9). Numerical Schlieren images, density and pressure cuts at time $t = 106$ μs, 204 μs, 301 μs. The red line is the iso-zero of the level set function. (For interpretation of the references to color in this figure legend, the reader is referred to the web version of this article.)

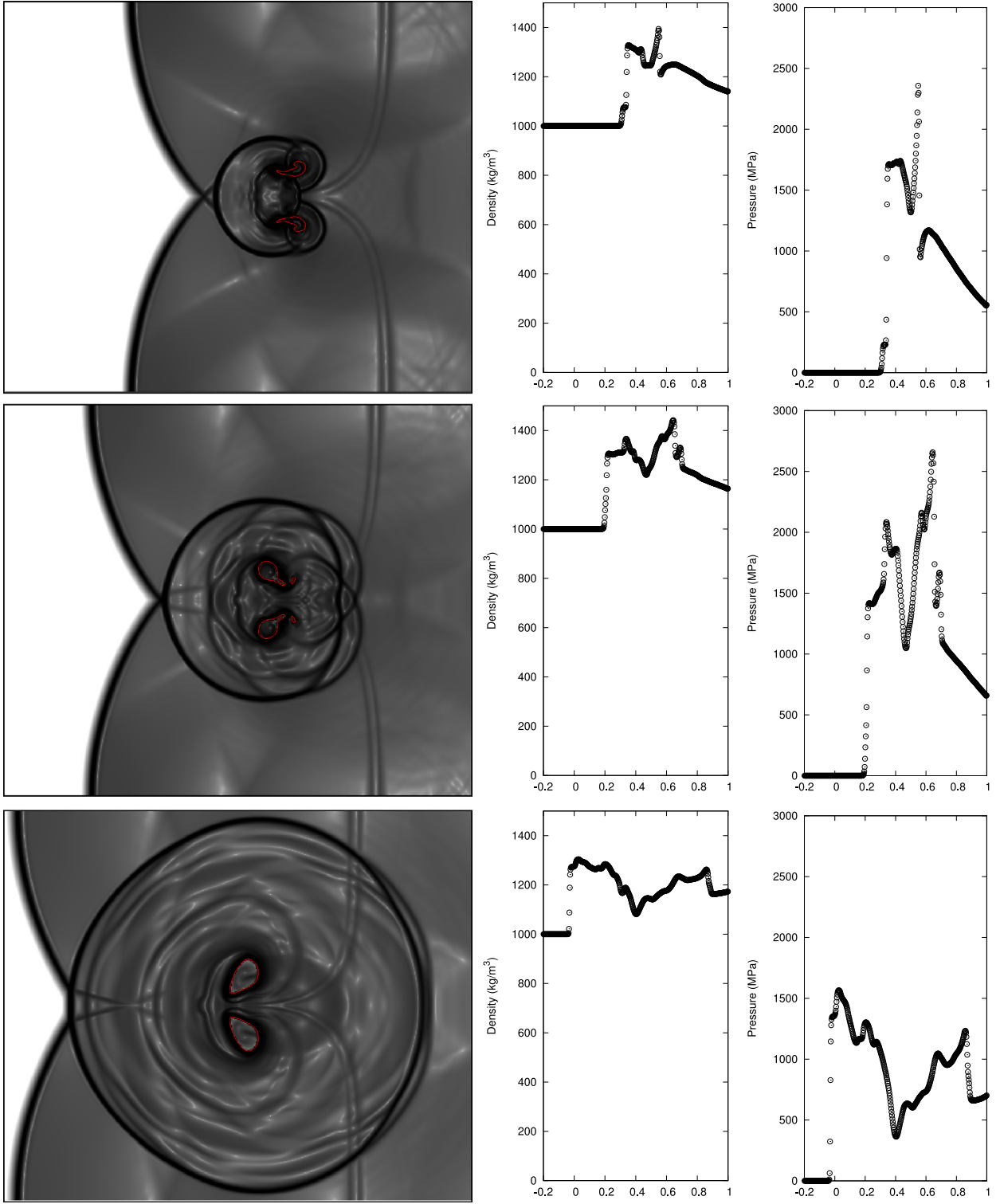

**Fig. 20.** Interaction of a Ma 1.422 water shock and an air bubble (TC9). Numerical Schlieren images, density and pressure cuts at time $t = 358$ μs, 406 μs, 500 μs. The red line is the iso-zero of the level set function. (For interpretation of the references to color in this figure legend, the reader is referred to the web version of this article.)

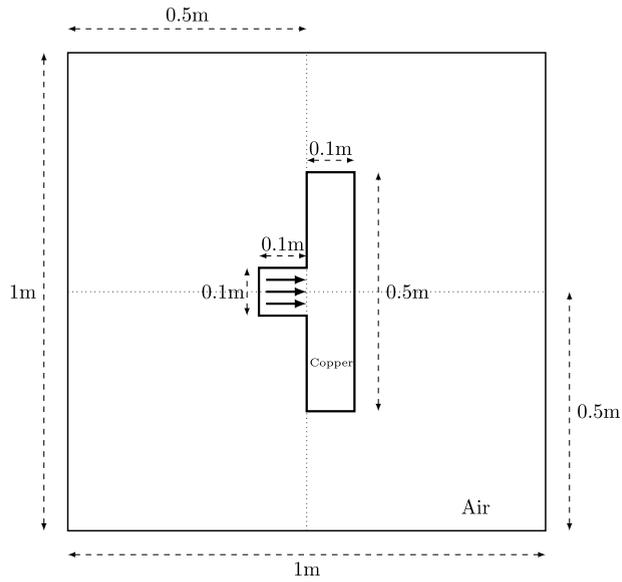

**Fig. 21.** Sketch of the initial configuration for the impact test case TC10 and TC11.

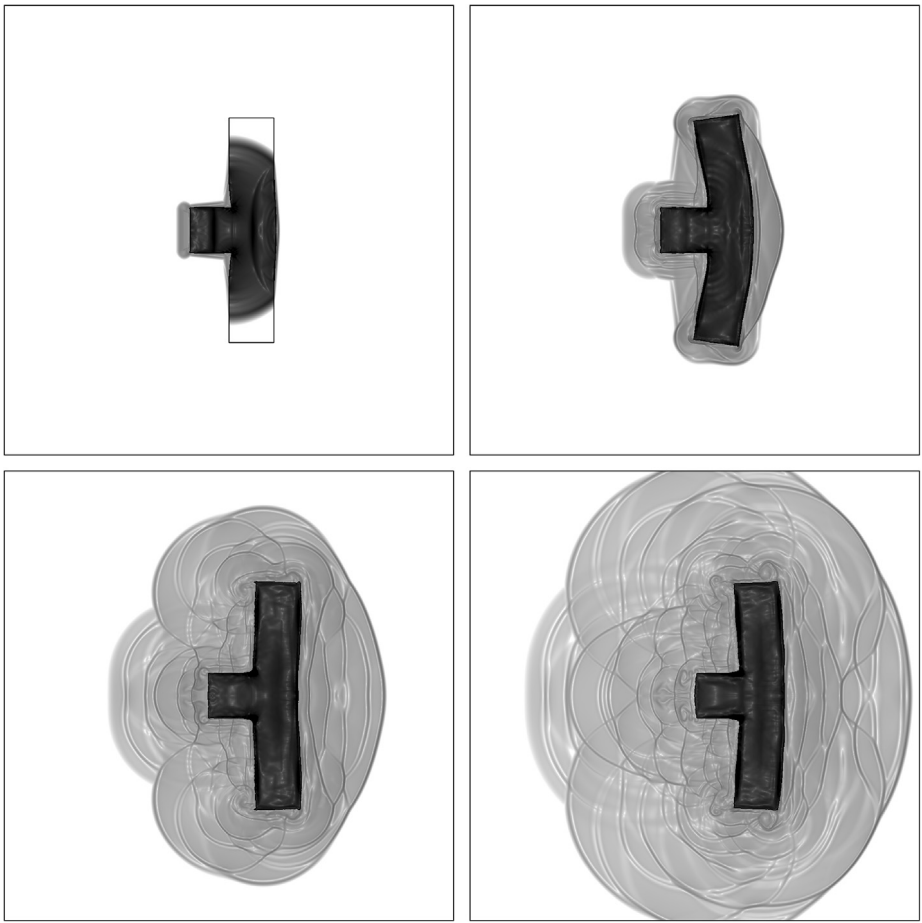

**Fig. 22.** Impact of a projectile on a plate for $\chi = 5 \cdot 10^{10}$ Pa (TC10). Numerical Schlieren pictures at $t = 27$ μs, 140 μs, 420 μs, 710 μs.

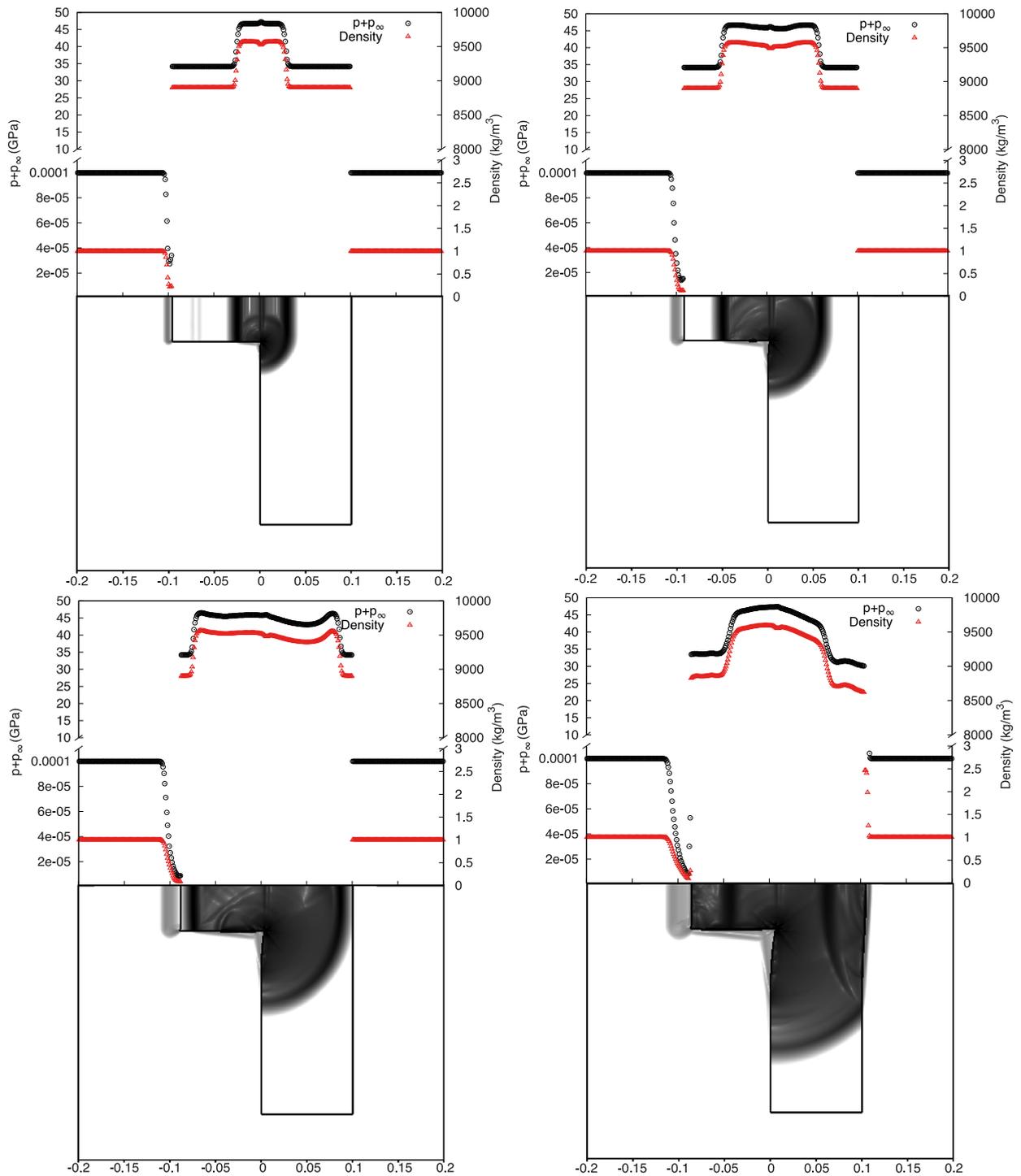

**Fig. 23.** Impact of a projectile on a plate for $\chi = 5 \cdot 10^{10}$ Pa (TC10). Numerical Schlieren pictures and section on the horizontal plane of symmetry for the pressure (black circles) and the density (red triangles) at $t_1 = 5$ µs, $t_2 = 10$ µs, $t_3 = 15$ µs, $t_4 = 25$ µs. (For interpretation of the references to color in this figure legend, the reader is referred to the web version of this article.)

In Fig. 22, one can see numerical Schlieren pictures of the flow for TC10. We can distinguish the propagation of a wave in the solid and its transmission in air. The elastic material is deformed and oscillates while being displaced rightward. This behavior is in good agreement with the literature, see [9,20]. In [9], the wave transmitted in air seems faster than it should compared to our results and along directions which seem to be influenced by the presence of a small amount of elastic solid mixed with air. Here, the two phases are pure and the waves travel without any perturbation.

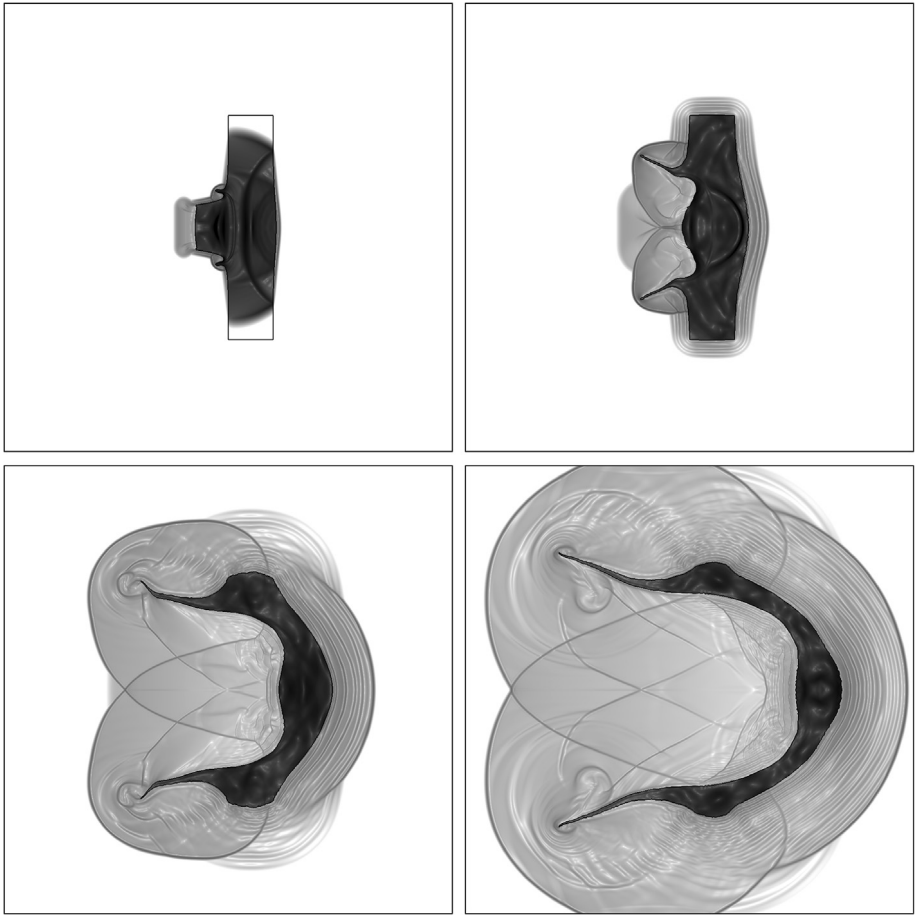

**Fig. 24.** Impact of a projectile on a plate for $\chi = 0$ (TC11). Numerical Schlieren pictures at $t = 38$ μs, 140 μs, 420 μs, 710 μs.

In Fig. 23 the density and the pressure on a section through the horizontal symmetry axis are shown with the corresponding Schlieren pictures. The abscissa $x = 0$ corresponds to the initial velocity discontinuity. At $t = 0$, for $-0.1 < x < 0.1$ the phase is elastic, elsewhere there is perfect gas. Two shock waves emanate left and rightward from the impact point ($t_1$ and $t_2$) and an expansion fan develops in the fluid to the left ($x < -0.1$). In snapshots $t_3$ and $t_4$ one can observe the reflection of the shock waves at the borders of copper as rarefaction waves, the expansion fan that continues to develop in the fluid for $x < -0.1$ and a shock wave transmitted in the fluid for $x > 0.1$. Overall, the density plot shows that the ensemble projectile/plate is displaced to the right.

In Fig. 24, one can see the results for $\chi = 0$ (TC11), this corresponds to the fluid limit. The plate undergoes extreme deformation as there is no force to bring back the structure to its reference configuration. The solution is qualitatively in good agreement with [9], although the two fluids are significantly mixed in that simulation.

### 8.2.2.1. Consistency between $\rho$ and $\nabla_x Y$

Taking $\xi = Y(x, t)$ in Eq. (5) we have

$$\rho(x, t) = \det\big(\nabla_x Y(x, t)\big) \rho_0\big(Y(x, t)\big) \tag{50}$$

therefore, if the initial density $\rho_0$ is constant, the equation of mass conservation is actually redundant with the equation on $\nabla_x Y$.

The quantity $|\rho/\rho_0 - \det(\nabla_x Y)|$ is a measure of the error between the continuity equation and the equation on $\nabla_x Y$. We show this field inside the solid at final time and its evolution from initial time in different norms in Figs. 25 and 26 for TC10. The errors are produced at geometrical singularities and are small elsewhere. Moreover, the error stabilizes in time.

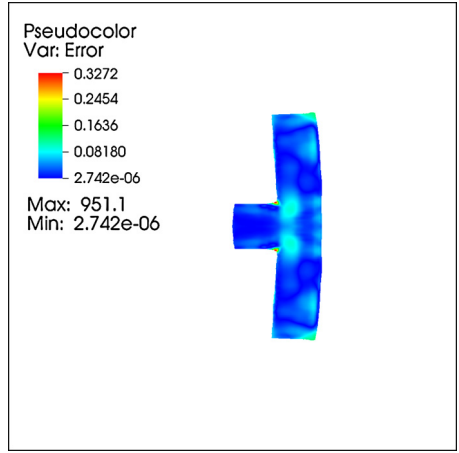

**Fig. 25.** TC10: $|\rho/\rho_0 - \det(\nabla_x Y)|$ field at final time (TC10).

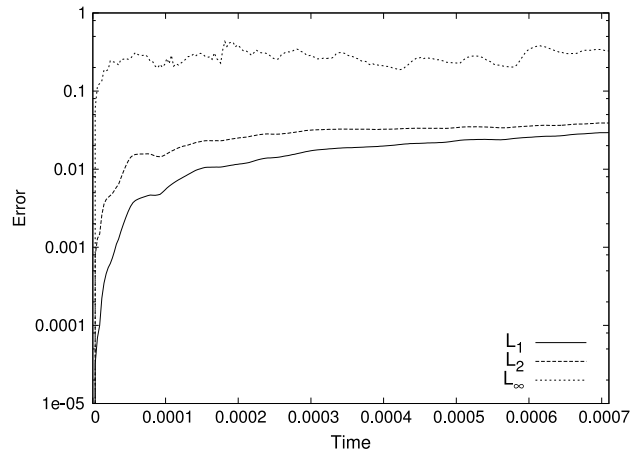

**Fig. 26.** Continuity error evolution in time for different norms (TC10).

## 9. Conclusions

We presented a simple method to deal with compressible multimaterials. This approach is simpler with respect to existing schemes as no ghost material is defined at the interface and no mixture model is needed.

This scheme is based on a fully Eulerian model that we have systematically presented. We have explicitly derived the characteristic speeds and the conditions of hyperbolicity for the neohookean constitutive law we use.

The computations in one and two space dimensions show that the contact discontinuity separating different materials remains sharp and that the numerical solution is in good agreement with the exact one or with numerical solutions presented in the literature, even when large velocity, pressure and material differences characterize the flow.

Extensions to three-dimensional configurations and to models including surface effects and plasticity are under current investigation.

### Acknowledgements

This study has been carried out with financial support from the French State, managed by the French National Research Agency (ANR) in the frame of the "Investments for the future" Programme IdEx Bordeaux (ANR-10-IDEX-03-02), Cluster of excellence CPU.